\documentclass{article}

\title{Derived autoequivalences from periodic algebras}
\author{Joseph Grant\\Graduate School of Mathematics,
Nagoya University,\\Chikusa-ku, Nagoya, 464-8602, Japan\\\texttt{joseph.grant@gmail.com}}
\date{}

\usepackage{latexsym}
\usepackage{verbatim}
\usepackage{amsmath}
\usepackage{amssymb}
\usepackage{mathrsfs}
\usepackage{amsthm}
\input xy
\xyoption{all}
\newcommand{\be}{\begin{equation}}
\newcommand{\ee}{\end{equation}}

\newcommand{\cma}{\textrm{,}}
\newcommand{\fs}{\textrm{.}}
\newcommand{\semic}{\textrm{;}}

\newcommand{\da}{\textrm{-}}
\newcommand{\dba}{{\Db(A)}}
\newcommand{\Hom}{\operatorname{Hom}\nolimits}
\newcommand{\RHom}{\operatorname{\textbf{R}Hom}\nolimits}
\newcommand{\dert}{\otimes^\textbf{L}}

\newcommand{\Aut}{\operatorname{Aut}\nolimits}
\newcommand{\Com}{\operatorname{Ch^b}\nolimits}

\newcommand{\cone}{\operatorname{cone}\nolimits}

\newcommand{\add}{\textbf{-add}}
\newcommand{\fa}{\textrm{ for all }}
\newcommand{\ev}{\textrm{ev}}

\newcommand{\Z}{\mathbb{Z}}

\newcommand{\N}{\mathbb{N}}
\newcommand{\BP}{\mathbb{P}}
\newcommand{\id}{\textrm{id}}

\newcommand{\mMod}{\textbf{-mod}}

\newcommand{\A}{\mathcal{A}}
\newcommand{\B}{\mathcal{B}}
\newcommand{\CC}{\mathcal{C}}

\newcommand{\T}{\mathcal{T}}
\newcommand{\gen}[1]{\langle#1\rangle}

\newcommand{\onto}{\twoheadrightarrow}
\newcommand{\into}{\hookrightarrow}

\newcommand{\env}{{E^\textrm{en}}}
\newcommand{\eop}{{E^\textrm{op}}}
\newcommand{\rsa}{\rightsquigarrow}
\newcommand{\arr}[1]{\stackrel{#1}{\to}}
\newcommand{\bs}{\backslash}

\newcommand{\Db}{\operatorname{D^b}\nolimits}
\newcommand{\Chb}{\operatorname{Ch^b}\nolimits}
\newcommand{\Dbsub}[1]{\operatorname{D^b_{\textit{$ #1 $}}}\nolimits}
\newcommand{\Kb}{\operatorname{K^b}\nolimits}
\newcommand{\K}{\operatorname{K}\nolimits}
\newcommand{\rad}{\operatorname{rad}\nolimits}
\newcommand{\soc}{\operatorname{soc}\nolimits}
\newcommand{\coker}{\operatorname{coker}\nolimits}
\newcommand{\en}{\textrm{en}}
\newcommand{\rsar}[1]{\stackrel{#1}{\rightsquigarrow}}
\newcommand{\CD}{\mathcal{D}}
\newcommand{\arrl}[1]{\stackrel{#1}{\longrightarrow}}
\newcommand{\im}{\operatorname{im}\nolimits}
\newcommand{\End}{\operatorname{End}\nolimits}
\newcommand{\op}{{\operatorname{op}\nolimits}}

\def\octax#1#2#3#4#5#6#7#8#9{%
  \def\ArgI{{#1}}%
  \def\ArgII{{#2}}%
  \def\ArgIII{{#3}}%
  \def\ArgIV{{#4}}%
  \def\ArgV{{#5}}%
  \def\ArgVI{{#6}}%
  \def\ArgVII{{#7}}%
  \def\ArgVIII{{#8}}%
  \def\ArgIX{{#9}}%
  \octaxxptii
}
\def\octaxxptii#1#2#3#4#5#6#7#8#9{%
  \xymatrix@=4pt{
{}\phantom{ABCDE} & {}\phantom{.} & {}\phantom{ABCDE} & {}\phantom{.} & {}\phantom{ABCDE} & {}\phantom{.} & {}\phantom{ABCDE} & {}\phantom{.} & {}\phantom{ABCDE} & {}\phantom{.} & {}\phantom{.} \\
 & & & & & & & & & & & & {}\phantom{.} \\
 & & & & & & & & & & & & {}\phantom{.} \\
 & & & &{\ArgIII}\ar@{-->}[rrddd]^{#7}\ar@{~>}[uuurr]^{\ArgIX}& & & & & & & & {}\phantom{.} \\
 & & & & & & & & & & & & {}\phantom{.} \\
 & & & & & & & & & & & & {}\phantom{.} \\
 & &{\ArgII}\ar[uuurr]^{\ArgVIII}\ar[rrd]^{#4}& & & &{\ArgV}\ar@{-->}[rrddd]^{#8}\ar@{~>}[rrrruu]^{#3}& & & & & & {}\phantom{.} \\
 & & & &{\ArgIV}\ar[urr]^{#2}\ar[rrrrdd]^{#5}& & & & & & & & {}\phantom{.} \\
 & & & & & & & & & & & & {}\phantom{.} \\
{\ArgI}\ar[uuurr]^{\ArgVII}\ar[uurrrr]^{#1}& & & & & & & &{\ArgVI}\ar@{~>}[rrrrdd]^{#6}\ar@{~>}[rrddd]^{#9}& & & & {}\phantom{.} \\
 & & & & & & & & & & & & {}\phantom{.} \\
 & & & & & & & & & & & & {}\phantom{.} \\
 & & & & & & & & & & & & {}\phantom{.}   
}
}

\newtheorem{thm}{Theorem}[subsection]
\newtheorem{lem}[thm]{Lemma}
\newtheorem{prop}[thm]{Proposition}

\newtheorem{cor}[thm]{Corollary}
{\theoremstyle{definition} \newtheorem{defn}[thm]{Definition}
\newtheorem{rmk}[thm]{Remark}}
\newenvironment{pf}{Proof: \em}{ \hfill $\Box$}

\newtheorem*{thma}{Theorem A}
\newtheorem*{thmb}{Theorem B}

\begin{document}

\maketitle

\begin{abstract}
We present a construction of autoequivalences of derived categories of symmetric algebras based on projective modules with periodic endomorphism algebras.  This construction generalises autoequivalences previously constructed by Rouquier-Zimmermann and is related to the autoequivalences of Seidel-Thomas and Huybrechts-Thomas.  
We show that compositions and inverses of these equivalences are controlled by the resolutions of our endomorphism algebra and that each autoequivalence can be obtained by certain compositions of derived equivalences between algebras which are in general not Morita equivalent.

\emph{Keywords:} symmetric algebra,  derived equivalence, two-sided tilting complex, spherical twist, periodic algebra, derived Picard group, tilting mutation

\emph{2010 Mathematics Subject Classification:} 18E30, 16E35, 16D50
\end{abstract}

\tableofcontents
\section{Introduction}

Derived categories are increasingly important in many areas of mathematics including representation theory and algebraic geometry.  They capture much of the homological information of an abelian category.  For example, the Grothendieck group of an abelian category is preserved under derived equivalence, and derived equivalent algebras have the same Hochschild cohomology.  Autoequivalences of derived categories are symmetries of these categories, and so give us extra homological information about the objects we are studying.  They are important in the representation theory of finite groups, in the study of knot invariants, and in many situations where abelian and triangulated categorification plays a role.

Important explicit autoequivalences of derived categories were constructed by Rouquier and Zimmermann \cite{rz} in the representation theory of finite groups, and by Seidel and Thomas \cite{st}, under the name \emph{spherical twists}, in algebraic geometry.  These autoequivalences have a similar form: they both depend only on the existence of objects in the derived category which satisfy a certain Calabi-Yau property and have particular endomorphism or extension algebras.  Another example of derived autoequivalences, called \emph{$\BP^n$-twists}, was given by Huybrechts and Thomas \cite{ht}.  Again, this construction only requires a Calabi-Yau object with a suitable extension algebra and gives an explicit functor which is an autoequivalence.  The $\BP^n$-twists and spherical twists are connected: when $n=1$, the $\BP^1$-twist is obtained by applying a spherical twist twice.

We can give analogous versions of these results in the setting of representations of symmetric algebras.  Let $A$ be a symmetric algebra over an algebraically closed field $k$ and $P$ be a projective left $A$-module.  If $\End_A(P)\cong k[x]/\gen{x^2}$, then the complex of functors
$$P\otimes_k\Hom_A(P,-)\to -\cma$$
where the nonzero differential is given by the obvious evaluation map, gives an autoequivalence of the derived category of finitely-generated $A$-modules which preserves the triangulated structure.  If $\End_A(P)\cong k[x]/\gen{x^{n+1}}$ for $n\geq1$, then the complex of functors
$$P\otimes_k\Hom_A(P,-)\to P\otimes_k\Hom_A(P,-)\to -\cma$$
where the first nonzero differential is the map $x\otimes 1-1\otimes\Hom_A(x,-)$ and the second is again the evaluation map, again gives an autoequivalence, and in the case $n=1$ the second functor is isomorphic to the square of the first.  These statements will both follow as special cases of Theorem A.

The present work is an attempt to fit the autoequivalences discussed above into a more general theory and to give an explanation of the connection between the two functors, based on \cite{jthesis}.  Our main theorem is as follows:
\begin{thma}
Suppose $A$ is a finite dimensional symmetric $k$-algebra and $P$ is a projective $A$-module.  Then if the algebra $E=\End_A(P)^\op$ is periodic of period $n$, i.e., $\Omega^n_{E\otimes_k\eop}(E)\cong E$ for some $n\in\N$, then there is a derived autoequivalence $\Psi_P:\Db(A)\to\Db(A)$, which we call a \emph{periodic twist}, which acts as the $n$-fold suspension functor $[n]$ on $P$ and acts as the identity on all objects of $\Db(A)$ which are right orthogonal to $P$.
\end{thma}
This theorem will follow from Theorem \ref{mainthm} and Corollary \ref{pn}, and will be proved in Subsection \ref{pfofeq}.  In fact we do more: we explicitly describe the autoequivalence in terms of a two-sided tilting complex, and we show that it suffices to require the weaker condition that $E$ is twisted-periodic, i.e., it is periodic up to some automorphism of $A$.  We also show how these periodic twists, and their endomorphism algebras, are controlled by bimodule resolutions of our endomorphism algebra $E$.

The periodic twists come from a ``circle of equivalences'':
\begin{thmb}
For every periodic twist $\Psi_P$, there exist 
non-trivial derived equivalences $\Db(A^{(i)})\to \Db(A^{(i+1)})$ of algebras $A=A^{(0)}$, $A^{(1)}$, \ldots, $A^{(n)}=A$ for $i=0,1,\ldots,n-1$ such that the composition of these equivalences agrees with $\Psi_P$ on objects.
\end{thmb}
$$
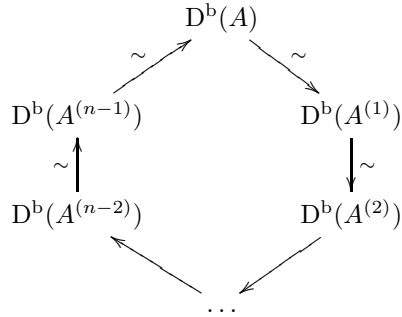
\xymatrix  @R=20pt @C=8pt  {
   & \Db(A)\ar[dr]^{\sim} & \\
 \Db(A^{(n-1)})\ar[ur]^{\sim} & & \Db(A^{(1)})\ar[d]^{\sim}\\
 \Db(A^{(n-2)})\ar[u]^{\sim} & & \Db(A^{(2)})\ar[dl] \\
 & {}\phantom{n}\:\hdots\:\phantom{n}\ar[ul] & 
}$$
This is Theorem \ref{circles}.  These derived equivalences are examples of tilting mutation (see, for example, \cite{ai}), as explained in Subsection \ref{combsec}.  While it is well-known that, in our situation, we can mutate our algebra $n$ times and obtain derived equivalences, Theorem B, as well as helping us to understand periodic twists, gives a condition for the series of mutations to result in an autoequivalence.

The paper is organised as follows.  In \textsection 2 we describe our notation and conventions, and review some background results that we will use later.  \textsection 3 describes the construction of the periodic twists in terms of tensoring with a two-sided tilting complex of $A\da A$-bimodules.  We calculate the action of our functor on the spanning class consisting of $P$ and elements of the derived category right orthogonal to it and use this to show that our functor is indeed an autoequivalence.  We also show that periodic twists are perverse, in the sense of Chuang and Rouquier.  In \textsection 4 we give explicit descriptions of the compositions and inverses of periodic twists, and prove some naturality results.  \textsection 5 contains a review of some basic tilting theory and a class of tilting complexes we call \emph{combinatorial}, and we use such complexes to construct our circle of equivalences and make the connection to periodic twists.  Finally, \textsection 6 contains some examples of periodic algebras and periodic twists.  In particular, Subsection \ref{truncex} examines the connections to the autoequivalences of \cite{rz}, \cite{st}, and \cite{ht} discussed above.

\textbf{Acknowledgements:}
I wish to thank Joe Chuang, who supervised me during my PhD studies, proposed this research project, and provided guidance and many suggestions during its development. I also wish to thank Jeremy Rickard and Osamu Iyama for suggestions and useful conversations, Alastair King for convincing me to write Section 4, Marcus Linckelmann for a suggestion that simplified the definition of the main functor presented in this article, and the anonymous referee for helpful comments and suggestions.

Most of the research presented here was conducted while the author was a PhD student at the University of Bristol, funded by the EPSRC, and Section 4 was completed, and the article was rewritten, while the author was a postdoctoral fellow at Nagoya University, funded by the JSPS.

\section{Background}

Let $k$ be an algebraically closed field and let $A$ be a finite-dimensional $k$-algebra.  All modules are left modules unless we specify otherwise, and the composition of the maps $f:X\to Y$ and $g:Y\to Z$ is written $gf$.  We will consider the following categories: the (abelian) category $A\mMod$ of finitely generated left $A$-modules, the (abelian) category of bounded (co)chain complexes of $A$-modules, denoted $\Com(A)$, the (triangulated) bounded homotopy category $\Kb(\CC)$ of an additive category $\CC$, and the (triangulated) bounded derived category of $A\mMod$, which we denote $\Db(A)$.  We will also consider the (additive) category $M\add$: this is the full subcategory with objects all objects isomorphic to summands of direct sums of $M$ if $M$ is an object of an abelian or triangulated category.

\subsection{Homological algebra}

We will use cohomological notation, so a complex $X=\bigoplus_{i\in\Z}X^i$ of $A$-modules will be of the form
$$\ldots\to X^{-2}\arr{d^{-2}}X^{-1}\arr{d^{-1}}X^0\arr{d^0}X^1\arr{d^1}X^2\to\ldots$$
and we say that $X^i$ is in degree $i$.    It will sometimes be convenient to write a complex with subscripts instead of superscrips, so $X_i=X^{-i}$, such as when writing projective resolutions of a module.  If we only display a finite section of a complex
$$L\to M\to\ldots\to N\cma$$
we assume that $N$ is in degree $0$ and that all modules to the left of $L$ or the right of $N$ are the zero module.

By $X[1]$ we denote the complex $X$ shifted to the left by one place with the sign of all differentials changed, i.e., 
$X[1]$ is the complex
$$\ldots\to X^{-1}\arr{-d^{-1}}X^0\arr{-d^0}X^1\arr{-d^1}X^2\arr{-d^{2}}X^3\to\ldots$$
with $X^1$ in degree $0$, and we define $X[n]$ for $n\in\Z$ following this convention.  If $M$ is an $A$-module, we will consider $M$ as a complex concentrated in degree $0$, and so for $n\in\Z$, $M[n]$ will denote the complex with $M$ in degree $-n$ and the zero module in every other degree.

In a triangulated category we will denote distinguished triangles by
$$X\arr{f} Y\arr{g} Z\stackrel{h}{\rsa}$$
which will mean we have morphisms as follows
$$X\arr{f} Y\arr{g} Z\arr{h} X[1]\fs$$
Recall that by the axioms of a triangulated category, given a triangle as above, we also have a shifted triangle
$$Y\arr{g} Z\arr{h} X[1]\stackrel{-f[1]}{\rsa}\fs$$
Also, given any morphism $f:X\to Y$ in $\Com(A)$, we have an object $\cone(f)$ and a short exact sequence
$$0\to Y\to \cone(f)\to X[1]\to0$$
in $\Com(A)$
which give rise to a distinguished triangle
$$X\arr{f} Y\to \cone(f)\rsa$$
in $\Db(A)$. 
For more information on triangulated and derived categories, see Chapter 10 of \cite{w}.

\subsection{Bimodules and symmetric algebras}

Given another $k$-algebra $B$, an $A\da B$-bimodule $M$ is a $k$-vector space which is simultaneously a left $A$-module and a right $B$-module such that the actions of $A$ and $B$ commute, i.e., for $a\in A$, $b\in B$, and $m\in M$, $(am)b=a(mb)$.  We will require that $k$ acts centrally on all $A\da B$-bimodules.  Note that $A$-modules can be identified with $A\da k$-bimodules, where $k$ denotes the trivial $k$-algebra.  We denote the category of $A\da B$-bimodules by $A\mMod\da B$, and the category of complexes of $A\da B$-bimodules by $\Com(A\da B)$.

The \emph{opposite algebra of $A$}, denoted $A^\op$, is the $k$-algebra with the same underlying vector space as $A$ but with multiplication defined by
$a\circ a'=a'a$
for $a,a'\in A$.  Right $A$-modules can be identified with $A^\textrm{op}$-modules, and $A\da B$-bimodules can be identified with left modules for the $k$-algebra $A\otimes_kB^\op$; we will use these identifications without further comment.  In the case $B=A$ we will call $A^\textrm{en}=A\otimes_kA^\op$ the \emph{enveloping algebra} of $A$.  Note that $A$ is an $A^\textrm{en}$-module, where the left and right action of $A$ are just given by the multiplication of $A$.

If $M$ is an $A$-module then we denote the \emph{$k$-linear dual} $\Hom_k(M,k)$ by $M^*$.  $M^*$ is a right $A$-module, with the right action of $A$ induced by the left action of $A$ on $M$, and as our modules are finite dimensional, $(M^*)^*$ and $M$ are naturally isomorphic as $A$-modules, so this is indeed a duality.  If $M$ is an $A\da B$-bimodule then $M^*$ is a $B\da A$-bimodule, and in particular, as $A$ is an $A\da A$-bimodule, so is $A^*$, with the left/right action on $A^*$ induced by the right/left action on $A$, respectively.

We say that $A$ is \emph{symmetric} if we have an isomorphism $A\cong A^*$ of $A\da A$-bimodules.  Hereafter we will assume that $A$ is symmetric.  Equivalent definitions of symmetric algebras can be found in Section 3 of \cite{r}, and a discussion of the properties of symmetric algebras can be found in Section 1.6 of \cite{b}.  In particular, note that symmetric algebras are self-injective, so the projective and the injective modules for $A$ coincide, and the head of any indecomposable projective $A$-module is isomorphic to its socle, i.e., $P/\rad P\cong \soc P$.  Symmetric algebras are closed under Morita equivalence and derived equivalence \cite[Corollary 5.3]{rictwo}.  Also, if $A$ is symmetric, the functors $(-)^*=\Hom_k(-,k)$ and $(-)^\vee=\Hom_A(-,A)$ are isomorphic (see, e.g., Theorem 3.1 of \cite{r}). Examples of symmetric algebras include group algebras of finite groups, trivial extension algebras of finite-dimensional algebras, and cohomology algebras of compact oriented manifolds.

\subsection{Twisted modules}\label{twistedmod}

Given an $A$-module $M$ and an automorphism $\sigma$ of the algebra $A$, we define a \emph{twisted module} ${_\sigma M}$ where the action of $A$ is given by $a\circ m=\sigma(a)m$ for $a\in A$ and $m\in M$.  Similarly, if $M$ is an $A\da A$-bimodule and $\sigma,\tau\in\Aut(A)$, we define the twisted module ${_\sigma M_\tau}$ with action $a\circ m\circ a'=\sigma(a)m\tau(a')$.  $M_\sigma$ denotes the twisted bimodule ${_\id M_\sigma}$.  Note that for any $\sigma\in\Aut(A)$, twisting by $\sigma$ does not change the isomorphism class of the regular $A$-module, as $\sigma$ induces an $A$-module automorphism $A\to{_\sigma A}$.  However, if ${_\sigma A}$ denotes the $A\da A$-module $A$ with left action twisted by $\sigma$ and with the usual right action, then $A$ and ${_\sigma A}$ are isomorphic as $A\da A$ bimodules if and only if $\sigma$ is an inner automorphism.

For any $A$-module $M$, we have an isomorphism ${_\sigma M}\cong {_\sigma A}\otimes_A M$ of left $A$-modules.  This comes from a natural isomorphism
$$\tau^{A,\sigma}:{_\sigma(-)}\arr{\sim}{_\sigma A}\otimes_A-$$
of functors, which, when $\sigma=1$, specialises to the natural isomorphism
$$\tau^{A}:\id_A\arr{\sim}A\otimes_A-\fs$$

\subsection{Adjunctions}\label{adjs}
We refer to Chapter 4 of \cite{mac} for general information on adjunctions.  We write $F\dashv G$ if $F$ is left adjoint to $G$.

If $A$, $B$, and $C$ are finite-dimensional $k$-algebras and
$L$ is a $C\da A$-bimodule,
$M$ is a $A\da B$-bimodule, and
$N$ is a $C\da B$-bimodule,
we have the usual tensor-hom adjunctions giving natural isomorphisms
$$\Hom_{C\otimes_kB^\op}(L\otimes_AM,N)\cong\Hom_{A\otimes_k B^\op}(M,\Hom_C(L,N))$$
and
$$\Hom_{C\otimes_kB^\op}(L\otimes_AM,N)\cong\Hom_{C\otimes_k A^\op}(L,\Hom_{B^\op}(M,N))$$
of vector spaces.  If we consider only homomorphisms of left modules then we have an isomorphism of $B\da B$-bimodules
$$\Hom_{C}(L\otimes_AM,N)\cong\Hom_{A}(M,\Hom_C(L,N))$$
and an isomorphism of $C\da C$-bimodules
$$\Hom_{B^\op}(L\otimes_AM,N)\cong\Hom_{A^\op}(L,\Hom_{B^\op}(M,N))\fs$$

These adjunctions are natural not only in $M$ and $N$, but also in $L$: if $h:K\to L$ is a map of $C\da A$-bimodules, then the following diagram commutes:
$$\xymatrix{ \Hom_{C}(L\otimes_AM,N)\ar[r]^(.45){\sim}\ar[d]^{\Hom_{C}(h\otimes_AM,N)} & \Hom_{A}(M,\Hom_C(L,N))\ar[d]^{\Hom_{A}(M,\Hom_C(h,N))} \\
 \Hom_{C}(K\otimes_AM,N)\ar[r]^(.45){\sim} & \Hom_{A}(M,\Hom_C(K,N))
}$$

Similar adjuntions hold on the level of complexes and derived categories: for example,
$$\Hom_{\Db(C\otimes_kB^\op)}(L\dert_AM,N)\cong\Hom_{\Db(A\otimes_k B^\op)}(M,\RHom_C(L,N))\fs$$

If $M$ and $P$ are $A$-modules and $P$ is projective, then
$$P^\vee\otimes_AM\cong\Hom_A(P,M)$$
so in particular, as $A$ is symmetric,
$$P^*\otimes_AM\cong\Hom_A(P,M)$$
This implies that if $A$ and $B$ are symmetric and $X$ is a bounded complex of $A\da B$-modules which are projective both as left $A$-modules and as right $B$-modules, then $X\otimes_A-$ is both left and right adjoint to $X^*\otimes_A-$.  Note that due to our conditions on $X$, these tensor products of complexes are their own derived functors.

\subsection{Endomorphism algebras}

In what follows, we will consider a distinguished projective $A$-module, which we will denote $P$, and its endomorphism algebra $\End_A(P)$.  We write $E=\End(P)^\op$, and then $P$ is an $A\da E$-bimodule with the obvious right action of $E$.  If $P=\bigoplus_{i=1}^rP_i$, note that $E=\bigoplus_{i=1}^r \Hom_A(P,P_i)$, so each $\Hom_A(P,P_i)$ is a projective (left) $E$-module.  Then the functor
$$\Hom_A(P,-):P\add\to E\mMod$$
sends $P_i$ to $\Hom_A(P,P_i)$ and induces an equivalence of additive categories
$$P\add\arr{\sim}E\add\fs$$

\section{Constructing autoequivalences}
We describe the construction of a two-sided tilting complex based on objects with twisted-periodic endomorphism algebras.  These tilting complexes induce autoequivalences which we call \emph{periodic twists}.

\subsection{Periodic algebras}
Let $E$ be a $k$-algebra.

\begin{defn}\label{def-twper}
We say that the algebra $E$ is \emph{twisted-periodic}, or \emph{$\sigma$-periodic}, if
there exists an integer $n\geq1$, an automorphism $\sigma$ of the algebra $E$, and an exact sequence
$$P_{n-1}\arr{d_{n-1}}P_{n-2}\to\ldots\to P_2\arr{d_2}P_1\arr{d_1}P_0$$
of projective $\env$-modules such that $\coker d_1\cong E$ and $\ker d_{n-1}\cong {E_\sigma}$.
We call $n$ the \emph{period of $E$}
and call the complex
$$P_{n-1}\arr{d_{n-1}}P_{n-2}\to\ldots\to P_1\arr{d_1}P_0$$
the \emph{truncated resolution of $E$}, denoted $Y$.  If $\sigma$ is the identity automorphism then we say $E$ is \emph{periodic}.
\end{defn}

Note that we do not insist that period $n$ is minimal.  However, every $\sigma$-periodic algebra has a minimal period, which depends on $\sigma$.

We record the following theorem, which is very useful for testing whether a given algebra is twisted periodic.  The first statement of the theorem is due to M. C. R. Butler.
\begin{thm}[\cite{gss}]\label{simpleper}
All twisted periodic algebras are self-injective, and moreover an algebra $E$ is twisted periodic if and only if every simple $E$-module has a periodic projective resolution.
\end{thm}
More information on periodic algebras can be found in the survey by Erdmann and Skowro\'nski \cite{erdsko}.

The following result will be important in Section 4.
\begin{prop}\label{truncresisom}
If a twisted-periodic algebra has two truncated resolutions of the same length then they are chain homotopy equivalent.

\begin{pf}
Suppose we have truncated resolutions
$$P_{n-1}\arr{d_{n-1}}P_{n-2}\to\ldots\to P_1\arr{d_1}P_0$$
and
$$Q_{n-1}\arr{\partial_{n-1}}Q_{n-2}\to\ldots\to Q_1\arr{\partial_1}Q_0$$
of $E$, so $\coker d_1\cong\coker\partial_1\cong E$ and $\ker d_{n-1}\cong\ker\partial_{n-1}\cong E_\sigma$ for some $\sigma\in\Aut(E)$.  Then as both complexes are truncations of projective resolutions, by Theorem 1.1 of \cite{mannantrunc}, the complexes
$$P_{n-1}\oplus P'\stackrel{(d_{n-1},0)}{\longrightarrow}P_{n-2}\to\ldots\to P_1\arr{d_1}P_0$$
and
$$Q_{n-1}\oplus Q'\stackrel{(\partial_{n-1},0)}{\longrightarrow}Q_{n-2}\to\ldots\to Q_1\arr{\partial_1}Q_0$$
are chain homotopy equivalent, where $P'=Q_{n-1}\oplus P_{n-2}\oplus Q_{n-3}\oplus\ldots$ and $Q'=P_{n-1}\oplus Q_{n-2}\oplus P_{n-3}\oplus\ldots$.  But as $\ker d_{n-1}\cong\ker\partial_{n-1}$, repeated application of Schanuel's Lemma combined with the Krull-Schmidt Theorem shows that $P'\cong Q'$.
\end{pf}
\end{prop}

\subsection{Construction of the functor}\label{constr}
From now on, we will fix the following notation: let $A$ be a finite-dimensional symmetric algebra, $P$ be a projective $A$-module, and $E=\End_A(P)^\op$.

Suppose that $E$ is twisted periodic of period $n$.  Then $E$ has a truncated resolution $Y$, which we can consider as a chain complex of $E\da E$-bimodules with $Y_0$ concentrated in degree $0$.  We have a map $f:Y\to E$ of chain complexes of $E\da E$-bimodules given by the isomorphism $\coker d_1\cong E$.

Consider the following composition of maps in $\Chb(A\da A)$, which we denote $g$:
$$P\otimes_EY\otimes_EP^\vee\arr{P\otimes f\otimes P^\vee} P\otimes_EE\otimes_EP^\vee\arr{\sim} P\otimes_EP^\vee\arr{\ev}A$$
Here, the final map is given by the obvious evaluation map $p\otimes f\mapsto f(p)$.

Recall that we have a natural isomorphism $\tau^A:\id_A\arr{\sim}A\otimes_A-$ as described in subsection \ref{twistedmod}.  As $P$ is projective we have adjunctions $P\otimes_E-\dashv P^\vee\otimes_A-$ and $-\otimes_EP^\vee\dashv-\otimes_AP$.  We denote the counit of the first adjunction by $\varepsilon^L$ and the counit of the second by $\varepsilon^R$.  The next result shows that $g$ is unambiguously defined, and also provides a connection with the definition presented in \cite{jthesis}.
\begin{lem}
The following maps in $\Chb(A\da A)$ are equal:
\begin{enumerate}\label{altg2}
\item $\ev\circ (P\otimes_E(\tau^E_{P^\vee})^{-1})\circ (P\otimes_Ef\otimes_EP^\vee)$;
\item $\ev\circ ((\tau^{E^\op}_{P})^{-1}\otimes_EP^\vee)\circ (P\otimes_Ef\otimes_EP^\vee)$;
\item $\varepsilon^L_A\circ (P\otimes_E(\tau^{A^{\op}}_{P^\vee}\circ (\tau^E_{P^\vee})^{-1}))\circ (P\otimes_Ef\otimes_EP^\vee)$;
\item $\varepsilon^R_A\circ ((\tau^A_P\circ(\tau^{E^\op}_P)^{-1})\otimes_EP^\vee)\circ (P\otimes_Ef\otimes_EP^\vee)$;
\item\label{def5} the image of the map $f$ in the following chain of adjunctions: \begin{align*}
\Hom_{\Com(E\da E)}(Y,E) &\cong {} \Hom_{\Com(E\da E)}(Y,\Hom_A(P,P))\\
 & {} \cong \Hom_{\Com(A\da E)}(P\otimes_E Y,P)\\
 & {} \cong \Hom_{\Com(A\da E)}(P\otimes_E Y,\Hom_{A^\op}(P^\vee,A))\\
 & {} \cong \Hom_{\Com(A\da A)}(P\otimes_E Y\otimes_E P^\vee, A)
\end{align*}
\end{enumerate}

\begin{pf}
By naturality, the final map is given by the composition
$$P\otimes_EY\otimes_EP^\vee\to P\otimes_EE\otimes_EP^\vee\to P\otimes_E\Hom_A(P,P^{\vee\vee})\otimes_EP^\vee\arr{\ev\otimes P^\vee}P^{\vee\vee}\otimes_EP^\vee\arr{\ev}A\fs$$
Then all five maps begin with $P\otimes_E f\otimes_E P^\vee$, and we can check directly that the remaining compositions agree in the module category.
\end{pf}
\end{lem}

\begin{defn}
Let $X=\cone(g)\in\Com(A\da A)$.
\end{defn}

\begin{lem}\label{mostproj}
$P\otimes_EY\otimes_EP^\vee$ is a bounded complex of modules which are projective as $A\da A$-bimodules.

\begin{pf}
All components of $Y$ are projective $\env$-modules $Q$, so it suffices to show that $P\otimes_EQ\otimes_EP^\vee$ is projective as an $A^\textrm{en}$-module.  But as projective modules are just direct summands of free modules, it suffices to show that $P\otimes_E\env\otimes_EP^\vee=P\otimes_EE\otimes_kE^\op\otimes_EP^\vee\cong P\otimes_kP^\vee\cong P\otimes_kP^*$ is a projective $A\da A$-bimodule, which is clear.
\end{pf}
\end{lem}

\begin{cor}\label{xper}
$X$ is a bounded complex of modules which are projective as left and as right $A$-modules.

\begin{pf}
This follows as $X=\cone(P\otimes_EY\otimes_EP^\vee\arr{g}A)$, and the corresponding statement is true for $P\otimes_EY\otimes_EP^\vee$ by Lemma \ref{mostproj}, and is clearly true for the regular $A$-module.
\end{pf}
\end{cor}

As $X$ is a bounded complex whose terms are projective as right $A$-modules, the functor $X\otimes_A-:\Com(A)\to\Com(A)$ is exact, and so induces a functor on the derived category.
\begin{defn}
Given a symmetric algebra $A$ and a projective $A$-module $P$, if $E=\End_A(P)^\textrm{op}$ is periodic with truncated resolution $Y$ then the \emph{periodic twist at $P$}, denoted $\Psi_{P,Y}$, or $\Psi_P$ if the truncated resolution is understood, is the functor
$$\Psi_P=X\otimes_A-:\Db(A)\to\Db(A)\fs$$
\end{defn}

\begin{thm}\label{mainthm}
If $A$ is symmetric, $P$ is a projective $A$-module, and $E=\End_A(P)^\op$ is twisted periodic of period $n$, then the periodic twist $\Psi_P:\Db(A)\to \Db(A)$ is an equivalence of triangulated categories.
\end{thm}
We will prove this theorem in the next subsection.

\begin{rmk}
\begin{enumerate}
\item Note that periodic twists are not periodic, i.e., there is no positive integer $\ell$ such that $\Psi_P^\ell\cong\id$.  One way to see this is to use Corollary \ref{pn} below. 
\item The functor $\Psi_{P,Y}$ depends not only on $P$ and $Y$, but also on the choice of isomorphism $\coker d_1\cong E$.  However, different choices of isomorphism give isomorphic functors.  Moreover our map $\coker d_1\cong E$ will usually be clear from the context.
\end{enumerate}
\end{rmk}

Periodic twists only depend on the summands of $P$ and not on their multiplicity.
\begin{lem}\label{nomult}
Let $P=\bigoplus_{j\in J}P_j$ be a sum of non-isomorphic indecomposable projective $A$-modules and let $P'=\bigoplus_{j\in J}P_j^{m_j}$ with $m_j\geq1$ for all $j\in J$.  Then $E=\End_A(P)^\op$ is twisted-periodic if and only $E'=\End_A(P')^\op$ is twisted-periodic, and in this case the periodic twists $\Psi_P$ and $\Psi_{P'}$ are naturally isomorphic.

\begin{pf}
Let $M$ be the $E\da E'$-bimodule $P^\vee\otimes_AP'\cong\Hom_A(P,P')$.  On restricting to the left action of $E$, $M$ is clearly a projective generator of $E\mMod$ and so induces a Morita equivalence $E'\mMod\arr{\sim}E\mMod$.  By Morita theory, $M\otimes_{E'}M^\vee\cong E$ and $M^\vee\otimes_EM\cong E'$.  We have a map $\varepsilon_{P'}:P\otimes_EP^\vee\otimes_AP'\arr{}P'$, which is an isomorphism because $\varepsilon_P$ is an isomorphism and $P'\in P\add$, and so $P\otimes_EM\cong P'$ and therefore also $P\cong P'\otimes_{E'}M^\vee$.

If $E$ is twisted periodic, i.e., we have a short exact sequence $$0\to E_\sigma[n-1]\into Y\onto E\to0$$ in $\Chb(E\da E)$, then by applying the functor $M^\vee\otimes_E-\otimes_EM$ we get a short exact sequence isomorphic to $$0\to E'_{\sigma'}[n-1]\into Y'\onto E'\to0$$ in $\Chb(E'\da E')$, where $Y'= M^\vee\otimes_EY\otimes_EM$ and $\sigma'$ is induced by $\sigma$.  Hence if $E$ is twisted-periodic then so is $E'$.

By our definitions, $P'\otimes_{E'}Y'\otimes_{E'}(P')^\vee\cong P\otimes_EY\otimes_EP^\vee$ and we have a commuting square
$$\xymatrix{
P'\otimes_{E'}Y'\otimes_{E'}(P')^\vee\ar[r]^(.72){g'}\ar[d]^{\sim} & A\ar@{=}[d]\\
P\otimes_EY\otimes_EP^\vee\ar[r]^(.7)g & A
}$$
so the complexes of bimodules $X$ and $X'$ are isomorphic, and hence $\Psi_P$ and $\Psi_{P'}$ are naturally isomorphic.  An analogous argument works in the other direction.
\end{pf}
\end{lem}
Hence we will assume hereafter that $P=\bigoplus_{j\in J}P_j$ is a sum of pairwise non-isomorphic projective modules.

We may assume the algebra $A$ is indecomposable:
\begin{prop}\label{blocks}
Suppose the algebra $A$ can be decomposed into blocks as $A=A_1\times A_2\times\ldots\times A_c$ and write $P=P_{(1)}\oplus P_{(2)}\oplus\ldots \oplus P_{(c)}$ where $P_{(i)}$ lies in the block $A_i$, and write $E_i=\End_{A_i}(P_{(i)})$.  Then $E$ is twisted-periodic of period $n$ if and only if so is each $E_i$, and $\Psi_P\cong \Psi_{P_{(1)}}\times\ldots\times\Psi_{P_{(c)}}$.

\begin{pf}
Suppose $A=A_1\times A_2\times\ldots A_c$, and that $1_A=a_1+a_2+\ldots+a_c$ is the corresponding decomposition of the identity of $A$ into primitive orthogonal central idempotents.  It follows that $E=E_1\times E_2\times\ldots\times E_c$ and  $1_E=e_1+e_2+\ldots+e_c$ is the corresponding decomposition of the identity of $E$, where $e_i$ is the endomorphism of $P$ given by $p\mapsto a_ip$.  Suppose each $E_i$ has a periodic resolution as an $E_i\otimes_kE_i^\op$-module given by applying $e_i$ to the resolution of $E$ as an $E\otimes_kE^\op$-module.  Then each $E_i$ has a truncated resolution $Y_i$ with periodicity $n$ and so for each $i=1,\ldots,c$, by the construction in Section \ref{constr} we get a complex of $A\da A$-bimodules $X_i$.  Conversely, if each $E_i$ is twisted periodic of period $n$ we can sum the resolutions to show that $E$ is twisted periodic of period $n$ and we obtain a complex $X$.

We now show that $X\cong\bigoplus_{i=1}^cX_i$.  This follows as for all $i\in\{1,2,\ldots,c\}$, $P_{(i)}=P_{(i)}e_i$, $Y_i=e_iY_ie_i$, and $P_{(i)}^\vee=\Hom_A(a_iP_{(i)},A)=e_iP_{(i)}^\vee$.  So $$P_{(i)}\otimes_EY_j\otimes_EP_{(k)}^*=P_{(i)}e_i\otimes_Ee_jY_je_j\otimes_Ee_kP_{(k)}^\vee=P_{(i)}\otimes_Ee_ie_jY_je_je_k\otimes_EP_{(k)}^\vee$$
is nonzero if and only if $i=j=k$.  Finally we note that as $X\cong\bigoplus_{i=1}^c X_i$, we have $\bigoplus_{i=1}^cX_i\otimes_A-\cong X\otimes_A-\fs$
\end{pf}
\end{prop}

\subsection{Proof of the equivalence}\label{pfofeq}

Our proof that periodic twists are two-sided tilting complexes will follow Ploog's alternative proof that spherical twists are autoequivalences \cite[Theorem 1.27]{plo}.

A collection $S$ of objects in $\Db(A)$ is a \emph{spanning class} if for all $L\in\Db(A)$ the following two conditions hold:
\begin{enumerate}
\item if $\Hom_{\Db(A)}(K,L[i])=0$ for all $K\in S$ and all $i\in\Z$ then $L\cong0$;
\item if $\Hom_{\Db(A)}(L[i],K)=0$ for all $K\in S$ and all $i\in\Z$ then $L\cong0$.
\end{enumerate}
Given some $K\in\Db(A)$, its \emph{right orthogonal complement} is
$$K^\perp=\{L\in\Db(A) \; | \; \Hom_{\Db(A)}(K,L[i])=0\fa i\in\Z\}$$
and its \emph{left orthogonal complement} $^\perp K$ is defined similarly.  The following lemma can be found as Corollary 3.2 of \cite{r} and says in particular that projective modules for symmetric algebras are $0$-Calabi-Yau objects. 
\begin{lem}\label{perdual}
If $A$ is a finite-dimensional symmetric $k$-algebra and $P$ is a bounded complex of finitely generated projective $A$-modules, then for any $M\in\dba$, the vector spaces $\Hom_\dba(P,M)$ and $\Hom_\dba(M,P)$ are naturally dual.
\end{lem}
Hence if $P$ is a projective $A$-module then its left and right orthogonal complements coincide, so we will just refer to $P^\perp$ and call it the \emph{orthogonal complement}.

\begin{lem}$S=\{P\}\cup P^\perp$ is a spanning class for $\Db(A)$.

\begin{pf}Suppose we have some $L\in\dba$ such that $\Hom_{\Db(A)}(K,L[i])=0$ for all $K\in S$ and all $i\in\Z$.  Then choosing $K=P\in S$, we see that $\Hom_{\Db(A)}(P,L[i])=0$ for all $i\in\Z$, i.e., $L\in P^\perp\subset S$.  So by our starting assumption, choosing $K=L$ and $i=0$ gives $\id_L\in\Hom_\dba(L,L)=0$, so $L\cong0$.

Now suppose that $\Hom_{\Db(A)}(L[i],K)=0$ for all $K\in S$ and all $i\in\Z$.  Then choosing $K=P$ gives $\Hom_{\Db(A)}(L[i],P)=0$ for all $i\in\Z$.
By Lemma \ref{perdual}, $\Hom_{\Db(A)}(P,L[i])^*=\Hom_{D^b(A)}(L[i],P)$ so $\Hom_{\Db(A)}(P,L[i])=0$ for all $i\in\Z$.  So $L\in P^\perp$, and so by the same argument as before $L\cong0$.
\end{pf}
\end{lem}

We now investigate how the functor $\Psi_P$ acts on elements of our spanning class $S$.  It will be useful to consider the triangulated structure of our categories.  
As $E$ is twisted-periodic, we have a short exact sequence of chain complexes of $\env$-modules
$$0\to {E_\sigma}[n-1]\to Y\arr{f} E\to0$$
which gives rise to a distinguished triangle
$$Y\arr{f}E\to {E_\sigma}[n]\rsa$$
in $\Db(\env)$ which we will denote $\Delta$.
Also, as we defined $X$ as the cone of the map $g$ constructed in Subsection \ref{constr}, we have a distinguished triangle
$$P\otimes_E Y\otimes_E P^\vee\arr{g} A\to X\rsa$$
in $\Db(A^\textrm{en})$ which we will denote $\nabla$.

\begin{prop}\label{upispd}
The triangles $P\otimes_E\Delta$ and $\nabla\otimes_AP$ are isomorphic.

\begin{pf}
We first note that as $P$ is a projective left $A$-module, and as all modules in complexes of the triangle $\Delta$ are projective left $E$-modules, we do not need to derive the tensor products.

Applying $P\dert_E-$ to $\Delta$ gives
$$P\otimes_EY\arr{P\otimes_Ef} P\otimes_EE\to P\otimes_E{E_\sigma}[n]\rsa$$ and
applying $-\otimes_AP$ to $\nabla$ gives
$$P\otimes_EY\otimes_EP^\vee\otimes_AP \arr{g\otimes_AP} A\otimes_AP\to X\otimes_AP\rsa$$
so if we can find isomorphisms $\alpha,\beta$ such that the square
$$\xymatrix{
  P\otimes_EY\ar[r]^{P\otimes_Ef}\ar[d]^\alpha & P\otimes_EE\ar[d]^\beta \\
  P\otimes_EY\otimes_EP^\vee\otimes_AP\ar[r]^{\phantom{mmmnn}g\otimes_AP} & A\otimes_AP
}$$
commutes, then by the completion axiom and the 5-lemma for triangulated categories it will follow that $P\otimes_E\Delta\cong\nabla\otimes_AP$.

Let $\beta$ be the isomorphism $\tau^A_P\circ(\tau^{E^\op}_P)^{-1}:P\otimes_EE\arr{\sim}A\otimes_AP$ and let $\alpha=\eta^R_{P\otimes_EY}$, where $\eta^R$ is the counit of the adjunction $-\dert_EP^\vee\dashv-\otimes_AP$.  To see that $\alpha$ is an isomorphism, we first use the identity 
$$(\varepsilon^R_-\otimes_AP)\circ \eta^R_{-\otimes_AP}=\id_{-\otimes_AP}$$
applied to $P\otimes_EY\otimes_EP^\vee$, which tells us that $\eta^R_{P\otimes_EY\otimes_EP^\vee\otimes_AP}$ is a monomorphism.  By the naturality of $\eta^R$ we have a commuting diagram
$$\xymatrix  @R=30pt @C=80pt {
P\otimes_EY\otimes_EP^\vee\otimes_AP\ar[r]^(0.43){\eta^R_{P\otimes_EY\otimes_EP^\vee\otimes_AP}}\ar[d]^{\chi}  & P\otimes_EY\otimes_EP^\vee\otimes_AP\otimes_EP^\vee\otimes_AP\ar[d]^{\chi\otimes_EP^\vee\otimes_AP}\\
P\otimes_EY\ar[r]^(0.43){\eta^R_{P\otimes_EY}}  & P\otimes_EY\otimes_EP^\vee\otimes_AP
}$$
where $\chi$ is an isomorphism $P\otimes_EY\otimes_EP^\vee\otimes_AP\cong P\otimes_EY$ induced by $P^\vee\otimes_AP\cong E$.  As both vertical maps are isomorphisms, $\eta^R_{P\otimes_EY}$ must be a monomorphism.  As $P\otimes_EY$ is a perfect complex of $A\da E$-bimodules, we can work in the homotopy category.  Then as $P\otimes_EY\otimes_EP^\vee\otimes_AP\cong P\otimes_EY$, $\alpha$ induces a monic endomorphism of the complex of finite-dimensional modules $P\otimes_EY$, and hence $\alpha$ must be an isomorphism

It only remains to show that the above diagram commutes, i.e., that $$(g\otimes_AP)\circ\alpha = \beta\circ(P\otimes_Ef)\fs$$
As 
$g=\varepsilon^R_A\circ (\beta\dert_EP^\vee)\circ (P\dert_Ef\dert_EP^\vee)$
we have
$$(g\otimes_AP)\circ\alpha = \left((\varepsilon^R_A\circ (\beta\dert_EP^\vee)\circ (P\dert_Ef\dert_EP^\vee))\otimes_AP\right)\circ\eta^R_{P\otimes_EY}$$
which is equal to
$$(\varepsilon^R_A\otimes_AP)\circ ((\beta\dert_EP^\vee)\otimes_AP)\circ ((P\dert_Ef\dert_EP^\vee)\otimes_AP)\circ\eta^R_{P\otimes_EY}$$
by the functoriality of $-\otimes_AP$.  By the naturality of $\eta$, this is equal to
$$(\varepsilon^R_A\otimes_AP)\circ \eta^R_{A\otimes_AP}\circ \beta\circ (P\otimes_Ef)$$
and as $(\varepsilon^R_A\otimes_AP)\circ \eta^R_{A\otimes_AP}=\id_{A\otimes_AP}$ this is equal to $\beta\circ (P\otimes_Ef)$.
\end{pf}
\end{prop}

\begin{cor}\label{pn}
For all $K\in P^\perp$, $\Psi_P(K)\cong K$, and $\Psi_P(P)\cong P[n]$.

\begin{pf}
To calculate $X\otimes_AK$ for $K\in P^\perp$, we apply the functor $-\dert_AK$ to $\Delta$ to obtain the triangle
$$P\otimes_E Y\otimes_E P^\vee\otimes_AK\arr{g} K\to X\otimes_AK\rsa\fs$$
Consider the complex $P^\vee\otimes_AK\cong\Hom_A(P,K)$.  Note that $H^0(\Hom_A(P,K))\cong\Hom_{K^b(A)}(P,K)$, and as $P$ is a bounded complex of projective $A$-modules, $\Hom_{K^b(A)}(P,K)\cong\Hom_\dba(P,K)$, so $H^i(\Hom_A(P,K))=\Hom_\dba(P,K[i])$.  As $K\in P^\perp$, $\Hom_\dba(P,K[i])=0$ for all $i\in\Z$, so $\Hom_A(P,K)\cong P^\vee\otimes_AK$ has zero cohomology in every degree.  Therefore
$P^\vee\otimes_AK\cong0$ in $\Db(A)$ and so the map $K\to X\otimes_AK$ is an isomorphism.

To see that $X\otimes_AP\cong P[n]$, we just use Lemma \ref{upispd}: as $E_\sigma[n]$ is the third term in $\Delta$ and $X$ is the third term in $\nabla$, $X\otimes_AP\cong P\otimes_EE_\sigma[n]$.  We are only considering $E_\sigma[n]$ as a complex of left $E$-modules here, so it is isomorphic to $E[n]$, and so $X\otimes_AP\cong P\otimes_EE[n]\cong P[n]$.
\end{pf}
\end{cor}

Recall that we have an automorphism $\sigma$ of the algebra $E$ from the definition of a twiseted periodic resolution.  Suppose $P=\bigoplus_{j\in J}P_j^{m_j}$ with $m_j>0$ for all $j\in J$.  Then $\sigma$ will permute the set $\{e_j\}_{j\in J}$ of orthogonal idempotents $e_j\in E$ and we write $E_{\sigma(j)}=E\sigma(e_j)$.  Using the equivalence of categories $E\add\cong P\add$ we write $P_{\sigma(j)}$ for the summand of $P$ corresponding to $E_{\sigma(j)}$.

\begin{cor}\label{psigma}
$\Psi_P(P_j)=P_{\sigma(j)}[n]$.

\begin{pf}
From Proposition \ref{upispd} we see that
$$X\otimes_AP\cong P\otimes_E{E_\sigma }[n]\cong P_{\sigma}[n]\fs$$
Then $Pe_j=P_j^{m_i}$,
 so
$$X\otimes_APe_j\cong P_{\sigma}e_j[n]\cong (P\sigma(e_j))_{\sigma}[n]$$
and hence by the Krull-Schmidt property, $X\otimes_AP_j\cong (P_{\sigma(j)})_{\sigma}[n]$.  Finally, forgetting the right $E$-module structure, $X\otimes_AP_j\cong (P_{\sigma(j)})[n]$.
\end{pf}
\end{cor}

We say a $k$-linear triangulated category $\T$ can be \emph{decomposed} into triangulated categories $\T_1$, $\T_2\subset \T$ if $\T_1$ and $\T_2$ are nonzero, $\Hom_\T(\T_1,\T_2)=\Hom_\T(\T_2,\T_1)=0$, and any $U\in\T$ can be written as $V_1\oplus V_2$ with $V_i\in\T_i$.  If a triangulated category cannot be decomposed into two nonzero triangulated subcategories we call it \emph{indecomposable}.

Recall that an algebra can be decomposed as the direct product of its blocks.  We note the following easily proved result:
\begin{lem}\label{blockindec}
If $A=A_1\times A_2\times\ldots\times A_c$ is a decomposition of $A$ into blocks, then
\begin{enumerate}
\item\label{dec_block}  $\Db(A)$ can be decomposed into $\Db(A_1)$, $\Db(A_2)$, $\ldots$, $\Db(A_c)$;
\item\label{ind_block} if $c=1$, i.e., $A$ is an algebra with a single block, then $\Db(A)$ is indecomposable.
\end{enumerate}
\end{lem}

To apply the above calculations, we will need the following result:
\begin{thm}[Bridgeland]\label{bri}
If $F:\mathcal{C}\to\mathcal{D}$ is a functor between triangulated categories such that
\begin{enumerate}
 \item $\mathcal{C}$ is nonzero and $\mathcal{D}$ is indecomposable;
 \item there exists a functor $G:\mathcal{D}\to\mathcal{C}$ which is both left and right adjoint to $F$;
 \item for all objects $K,L$ in some spanning class for $\mathcal{C}$ and all $i\in\Z$ the natural homomorphisms
$$F:\Hom_\mathcal{C}(K,L[i])\to\Hom_\mathcal{D}(F(K),F(L[i]))$$
are bijective
\end{enumerate}
then $F$ is an equivalence of categories.

\begin{pf}
 This follows from Theorems 2.3 and 3.3 of \cite{br}.
\end{pf}
\end{thm}

{\emph{Proof of Theorem \ref{mainthm}}
By Proposition \ref{blocks} we can assume that $A$ has only one block.  Then the three conditions of Theorem \ref{bri} are satisfied as
\begin{enumerate}
 \item $\Db(A)$ is clearly nonzero and is indecomposable by Lemma \ref{blockindec};
 \item $X^*\otimes_A-$ is both left and right adjoint to $X\otimes_A-$ as $X$ is a bounded complex of left and right projective $A$-modules by Corollary \ref{xper};
 \item $X\otimes_A-$ is bijective on morphisms between elements of the spanning class $S=\{P\}\cup P^\perp$ by Corollary \ref{pn}
\end{enumerate}
so we have shown that $X\otimes_A-$ is an autoequivalence of $\Db(A)$. 
{ \hfill $\Box$}}

\subsection{Perverse equivalences}
We will show that periodic twists are perverse, as defined by Chuang and Rouquier \cite[Section 2.6]{rouqims}.

Let $A$ and $B$ be two derived equivalent algebras, and write $\A=A\mMod$ and $\B=B\mMod$.  Suppose that $I=\{S_i\}_{i=1}^r$ and $I'=\{S_i'\}_{i=1}^r$ are complete sets of isoclasses of simple modules for $A$ and $B$, respectively.  

If $J$ is a subset of $I$, then we have a Serre subcategory $\A_J$ of $\A$ generated by the simple modules $J$: this is the smallest full subcategory $\A_J$ of $\A$ containing the objects $J$ such that for each short exact sequence
$$0\to L\into M\onto N\to0$$
$M$ is an object of $\A_J$ if and only if so are $L$ and $N$.  In particular, $\A_\emptyset$ is the abelian category which contains only the zero object.  For such a Serre subcategory $\CC$, we denote by $\Dbsub{\CC}(A)$ the thick subcategory of $\Db(A)$ whose objects are complexes with homology contained in $\CC$.  We have a quotient $\Db(A)/\Dbsub{\CC}(A)$ which is defined by localizing $\Db(A)$ at all morphisms $f:U\to V$ such that for the distinguished triangle
$$U\arr{f}V\to W\rsa$$
we have $W\in\Dbsub{\CC}(A)$, and we have a functor $\Db(A)\to\Db(A)/\Dbsub{\CC}(A)$.  If $\CD$ is a Serre subcategory of $\B$ then an equivalence $\Dbsub{\CC}(A)\arr{\sim}\Dbsub{\CD}(B)$ induces an equivalence $\Db(A)/\Dbsub{\CC}(A)\arr{\sim}\Db(B)/\Dbsub{\CD}(B)$.

\begin{defn}[Chuang-Rouquier]
A derived equivalence $F:\Db(A)\to\Db(B)$ is \emph{perverse} if there exist filtrations
$$\emptyset=I_0\subset I_1\subset\ldots\subset I_t=I$$
and
$$\emptyset=I_0'\subset I_1'\subset\ldots\subset I_t'=I'$$
which define abelian subcategories $\A_i:=\A_{I_i}$ and $\B_i:=\B_{I'_i}$ of $\A$ and $\B$, and a function $p:\{1,\ldots,t\}\to\Z$, such that for each $1\leq i\leq t$,
\begin{enumerate}
 \item $F$ restricts to an equivalence of triangulated categories $\Dbsub{\A_i}(A)\arr{\sim}\Dbsub{\B_i}(B)$ and
 \item $F[-p(i)]$ induces a Morita equivalence $\A_i/\A_{i-1}\arr{\sim}\B_i/\B_{i-1}$, i.e., there exists an equivalence of categories making the following diagram of functors commute:
$$\xymatrix@=45pt{
 {}\Db(A)/\Dbsub{\A_{i-1}}(A)\ar[r]^\sim_{F[-p(i)]} & {}\Db(B)/\Dbsub{\B_{i-1}}(B) \\
 \A_i/\A_{i-1}\ar@{^{(}->}[u]\ar@{-->}[r]^{\exists\,\sim} & \B_i/\B_{i-1}\ar@{^{(}->}[u]
}$$
\end{enumerate}
\end{defn}

\begin{rmk}\label{pervdiag}
In order to understand the second condition of the above definition, it is useful to consider the following diagram, which can be visualised as a cube with two adjacent faces missing:
$$\xymatrix{
 &\Db(A)\ar[rr]^\sim_{F[-p(i)]}\ar@{->>}[dd] & &\Db(B)\ar@{->>}[dd] \\
 \A_i\ar@{^{(}->}[ur]\ar@{->>}[dd] & &\B_i\ar@{^{(}->}[ur]\ar@{->>}[dd] & \\
 & {}\Db(A)/\Dbsub{\A_{i-1}}(A)\ar'[r]^>>{\sim}[rr] & & {}\Db(B)/\Dbsub{\B_{i-1}}(B) \\
  \A_i/\A_{i-1}\ar@{^{(}->}[ur]\ar@{-->}[rr]^{\sim} & & \B_i/\B_{i-1}\ar@{^{(}->}[ur] &
}$$
The three squares which do not contain the dashed arrow always commute, and $F$ is perverse if we can define the dashed arrow in such a way that the bottom square also commutes.
\end{rmk}

Now suppose that $A$ is a symmetric algebra and $P$ is a projective $A$-module with twisted-periodic endomorphism algebra of period $n$.
\begin{prop}\label{pertwperv}
The periodic twist $\Psi_P:\Db(A)\arr{\sim}\Db(A)$ is perverse with both filtrations given by
$$\emptyset\subset I\bs J\subset I$$
where $J$ is the set of isoclasses of summands of $P/\rad P$ and the function $p$ is defined by $p(1)=0$ and $p(2)=n$.

\begin{pf}
First we show that if $V\in\Dbsub{\A_i}(A)$ then $\Psi_P(V)\in\Dbsub{\A_i}(A)$.  The only non-trivial case is $i=1$.  But $V\in\Dbsub{\A_i}(A)$ precisely when $\Hom_{\Db(A)}(P,V[m])=0$ for all $m\in\Z$, i.e., when $V\in P^\perp$, and so by Corollary \ref{pn}, $\Psi_P(V)=V$.  We also have $\Psi_P^{-1}(V)=V$, and so $\Psi_P$ restricts to an autoequivalence of $\Dbsub{\A_1}(A)$.

Next we show that $\Psi_P[-p(i)]$ induces autoequivalences of $\A_i/\A_{i-1}$ such that the above diagram commutes.  For $i=1$ we have $\Db(A)/\Dbsub{\A_{i-1}}(A)=\Db(A)$ and $\Psi_P$ acts as the identity on $\A_i$, so everything is clear.  For $i=2$ we consider an object $\bar{M}\in\A/\A_1$ and we will show that it is sent by the functors
$$\A/\A_i\into\Db(A)/\Dbsub{\A_1}(A)\arrl{\Psi_P[-n]}\Db(A)/\Dbsub{\A_1}(A)$$
to a complex with homology in only one degree, and therefore defines an object of $\A/\A_1$.

Choose some $A$-module $M$ which represents $\bar{M}$.  We will apply the composition of functors
$$ \A\stackrel{\iota}{\into}\Db(A)\arrl{\Psi_P[-n]}\Db(A)\stackrel{\pi}{\onto}\Dbsub{\A_1}(A)$$
where we have labelled the inclusion and projection by $\iota$ and $\pi$, respectively. 
Then the commutative diagram in Remark \ref{pervdiag} shows that applying these functors to $M$ is the same as applying the functors in the previous paragraph to $\bar(M)$.  We want to show that $\pi\circ\Psi_P\circ\iota(M)[-n]$ has homology concentrated in degree $0$, i.e., that the $i$th homology of $\Psi_P(M)[-n]$ is contained in $\A_i$ for $i\neq 0$.  This is true because, by Corollary \ref{pn},
\begin{align*}
 \Hom_{\Db(A)}(P,\Psi_P(M)[-n][m]) &\cong\Hom_{\Db(A)}(\Psi_P^{-1}(P)[n],(M)[m])\\
 &\cong \Hom_{\Db(A)}(P,M[m])
\end{align*}
which is zero if $m\neq0$.  So given any $\bar{M}\in \A/\A_1$, we can choose some $M$ representing $\bar{M}$ and send it to the preimage of $\pi\circ F(M)[-n]$ in the functor $$\A/\A_1\into\Db(A)/\Dbsub{\A_{1}}(A)\fs$$  This will define an endofunctor of $\A/\A_1$ making the necessary diagram commute, and by a similar procedure one can define an inverse, showing that this is an equivalence.
\end{pf}
\end{prop}

\section{Compatibility}

We show that the composition of functors $\Psi_P$ agrees with composition of truncated bimodule resolutions, and show that quasi-inverses have a similar interpretation.

\subsection{Multiple periods}
If we have two truncated resolutions $Y_1$ and $Y_2$ of a twisted-periodic algebra $E$ of periods $n_1$ and $n_2$, then we can splice them together to obtain a new truncated resolution $Y$ of $E$ as a twisted periodic algebra of period $n=n_1+n_2$.

The two truncated resolutions give us distinguished triangles $\Delta_1$
$$Y_1\arr{f_1}E\arr{e_1} E_{\sigma_1}[n_1]\rsar{d_1}$$
and $\Delta_2$
$$Y_2\arr{f_2}E\arr{e_2} E_{\sigma_2}[n_2]\rsar{d_2}$$
in $\Db(E^\en)$ which come from short exact sequences
$$0\to E_{\sigma_1}[n_1-1]\arr{-d_1[-1]}Y_1\arr{f_1}E\to0$$
and
$$0\to E_{\sigma_2}[n_2-1]\arr{-d_2[-1]}Y_2\arr{f_2}E\to0$$
of chain complexes of $E\da E$-bimodules, where we have used the same notation for maps in $\Com(E^\en)$ and $\Db(E^\en)$.  We apply the functor $(-)_{\sigma_1}[n_1-1]$ to the second short exact sequence and so get another short exact sequence
$$0\to E_{\sigma}[n-2]\arr{-(d_2)_{\sigma_1}[n_1-2]}(Y_2)_{\sigma_1}[n_1-1]\arr{(f_2)_{\sigma_1}[n_1-1]}E_{\sigma_1}[n_1-1]\to0$$
in $\Com(E^\en)$ where $\sigma=\sigma_2\sigma_1$ and $n=n_1+n_2$.  We obtain the following exact sequence in $\Com(E^\en)$:
$$\xymatrix@=12pt{
0\ar[rr] &&E_{\sigma}[n-2]\ar[rr]^{-(d_2)_{\sigma_1}[n-2]} && (Y_2)_{\sigma_1}[n_1-1]\ar[rd]_(0.4){(f_2)_{\sigma_1}[n_1-1]}\ar[rr] && Y_1\ar[rr]^{f_1} &&E\ar[rr] && 0\\
 &&&&{}\phantom{mmm}& E_{\sigma_1}[n_1-1]\ar[ru]_(0.6){-d_1[-1]} &{}\phantom{mmm}&&{}\phantom{mmm}&&
}$$
and put $Y=\cone((Y_2)_{\sigma_1}[n_1-1]\to Y_1)\in \Com(E^\en)$.

The truncated resolutions $Y_1$, $Y_2$, and $Y$ of $E$ give us three periodic twists $\Psi_{P,Y_1}$, $\Psi_{P,Y_2}$, and $\Psi_{P,Y}$.

\begin{prop}\label{comp}
There is a natural isomorphism $\Psi_{P,Y}\cong\Psi_{P,Y_2}\circ \Psi_{P,Y_1}$.

\begin{pf}
We break the proof into smaller steps.

\textit{Step 1:}
By rotating the distinguished triangle $\Delta_1$ in $\Db(E^\en)$ we obtain the distinguished triangle
\begin{equation}\label{qDel}
 E_{\sigma_1}[n_1-1]\arr{-d_1[-1]}Y\arr{f_1}E\rsar{e_1}\fs
\end{equation}
By applying the functor $-_{\sigma_1}[n_1-1]$ to $\Delta_2$ we get another distinguished triangle
\begin{equation}\label{DelSigShift}
(Y_2)_{\sigma_1}[n_1-1]\arr{(f_2)_{\sigma_1}[n_1-1]} E_{\sigma_1}[n_1-1]\arr{(e_2)_{\sigma_1}[n_1-1]}{E_{\sigma}}[n-1]\rsar{(d_2)_{\sigma_1}[n_1-1]}\fs 
\end{equation}
We defined $Y$ as $\cone(\delta:(Y_2)_{\sigma_1}[n_1-1]\to Y_1)$, where $\delta=-d_1\circ (f_2)_{\sigma_1}[n_1]$, so we have a distinguished triangle
\begin{equation}\label{defYtwo}
(Y_2)_{\sigma_1}[n_1-1]\arr{\delta}Y_1\arr{\iota}Y\rsar{\pi}
\end{equation}
coming from the short exact sequence of chain complexes
$$0\to Y_1\arr{\iota} Y\arr{\pi}(Y_2)_{\sigma_1}[n_1]\to0$$
where $\iota$ and $\pi$ are the obvious inclusion and projection.
Then by applying the octahedral axiom of the triangulated category $\Db(\env)$ to triangles (\ref{qDel}), (\ref{DelSigShift}), and (\ref{defYtwo})
as indicated in the following diagram:
$$\octax{(Y_2)_{\sigma_1}[n_1-1]}{E_{\sigma_1}[n_1-1]}{E_{\sigma}{[n-1]}}{Y_1}{Y}{E}{(f_2)_{\sigma_1}[n_1-1]}{(e_2)_{\sigma_1}[n_1-1]}{(d_2)_{\sigma_1}[n_1-1]}{}{\iota}{\pi}{-d_1[-1]}{f_1}{e_1}{}{f}{}$$
we obtain the distinguished triangle
$${E_{\sigma}}[n-1]\arr{}Y\arr{f}E\rsa\fs$$
We rotate this to obtain
$$Y\arr{f}E\arr{e}{E_{\sigma}}[n]\rsar{d}$$
which we denote $\Delta$.

\textit{Step 2:}
From the octahedral axiom we see that the composition
$$Y_1\arr{\iota}Y\arr{f}E$$
is equal to $f_1$, and so the map $f$ is just given by the cokernel of the last nonzero differential in the truncated resolution $Y$ of $E$.  So $f$ comes from a map of chain complexes, which we also call $f$, and we can use this map to construct
$$g:P\otimes_EY\otimes_EP^\vee\to A$$
as in Subsection \ref{constr} and then let $X=\cone(g)\in\Com(A\da A)$.  This gives a distinguished triangle
$$P\otimes_E Y\otimes_E P^\vee\arr{g} A\arr{h} X\rsar{i}$$
in $\Db(A\da A)$ which we will denote $\nabla$.

\textit{Step 3:}
We have a triangle $\nabla_1$ in $\Db(A^\en)$ where we label the maps as follows:
$$P\otimes_E Y_1\otimes_E P^\vee\arr{g_1} A\arr{h_1} X_1\rsar{i_1}$$
By rotating this triangle twice we obtain
\begin{equation}\label{rrNab}
X_1\arr{i_1}P\otimes_E Y_1\otimes_E P^\vee[1]\arr{-g_1[1]} A[1]\rsar{-h_1[1]}\fs
\end{equation}
We also obtain the triangle
\begin{equation}\label{PytwoP}
P\otimes_E(Y_2)_{\sigma_1}\otimes_EP^\vee[n_1]\arr{}P\otimes_EY_1\otimes_EP^\vee[1]\arr{}P\otimes_EY\otimes_EP^\vee[1]\rsar{}
\end{equation}
by applying $P\otimes_E-\otimes_EP^\vee$ to the triangle (\ref{defYtwo})
obtained in Step 1 and rotating three times.

We now apply $-\otimes_AX_1$ to $\nabla_2$ and define a distinguished triangle
\begin{equation}\label{altXDel}
P\otimes_E (Y_2)_{\sigma_1}\otimes_E P^\vee[n_1]\arr{\bar{g}} X_1\arr{\bar{h}} X_2\otimes_AX_1\rsar{\bar{i}} 
\end{equation}
using the following isomorphism of triangles:
$$\xymatrix{
   P\otimes_E (Y_2)_{\sigma_1}\otimes_E P^\vee[n_1]\ar[r]^{\phantom{spacing}\bar{g}}\ar[d]^{P\otimes_E\tau^{E^\op,\sigma_1}_{Y_2}\otimes_EP^\vee[n_1]} &  X_1\ar[r]^{\phantom{aa}\bar{h}\phantom{spac}}\ar[dd]^{\tau^{A}_{X_1}}  & X_2\otimes_AX_1\ar@{~>}[r]^{\phantom{aaa}\bar{i}}\ar@{=}[dd] & \\
   P\otimes_EY_2\otimes_EE_{\sigma_1}\otimes_EP^\vee[n_1]\ar[d]^{P\otimes Y_2\otimes \psi} & &\\
 P\otimes_EY_2\otimes_EP^\vee\otimes_AX_1\ar[r]^{\phantom{spacing}g_2\otimes X_1} & A\otimes_AX_1\ar[r]^{h_2\otimes X_1} & X_2\otimes_AX_1\ar@{~>}[r]^{\phantom{aaa}i_2\otimes X_1} &
}$$
Here, the isomorphism $\psi$ comes from the isomorphism 
$$\xymatrix{
 Y_1\otimes_EP^\vee\ar[r]^{f_1\otimes_EP^\vee}\ar[d]^{} & E\otimes_EP^\vee\ar[d]^{}\ar[r]^{e_1\otimes_EP^\vee} & E_{\sigma_1}\otimes_EP^\vee[n_1]\ar@{-->}[d]^{\psi}\ar@{~>}[r]^{\phantom{mmmmn}d_1\otimes_EP^\vee} & \\
  P^\vee\otimes_AP\otimes_EY_1\otimes_EP^\vee\ar[r]^{\phantom{mmmnn}P^\vee\otimes_Ag_1} & P^\vee\otimes_AA\ar[r]^{P^\vee\otimes_Ah_1} & P^\vee\otimes_AX_1\ar@{~>}[r]^{\phantom{mmm}P^\vee\otimes_Ai_1} &
}$$
of triangles $\Delta\otimes_EP^\vee\cong P^\vee\otimes_A\nabla$, which is proved similarly to Proposition \ref{upispd}.

Our aim is to apply the octahedral axiom of $\Db(A^\en)$ to the triangles (\ref{rrNab}), (\ref{PytwoP}), and (\ref{altXDel}) as indicated by the following diagram:
$$\octax{P(Y_2)_{\sigma_1} P^\vee[n_1]}{X_1}{X_2X_1}{PY_1P^\vee[1]}{PYP^\vee[1]}{A[1]}{\bar{g}}{\bar{h}}{\bar{i}}{-P\delta P^\vee[1]}{-P\iota P^\vee[1]}{-P\pi P^\vee[1]}{i_1}{-g_1[1]}{-h_1[1]}{}{}{}$$
where we have omitted the tensor product symbols to save space. 
To be able to do this, we need to check that
$$i_1\circ\bar{g}=-P\otimes_E\delta\otimes_EP^\vee\fs$$

\textit{Step 4:}
We check the above equality holds.  Recall that we defined $\delta=-d_1\circ (f_2)_{\sigma_1}[n_1]$ when constructing triangle (\ref{defYtwo}).

Consider
$$i_1\circ\bar{g}=i_1\circ(\tau^{A}_{X_1})^{-1}\circ(g_2\otimes_AX_1)\circ(P\otimes_EY_2\otimes_E\psi)\circ(P\otimes_E\tau^{E^\op,\sigma_1}_{Y_2}\otimes_EP^\vee[n])\fs$$
By naturality, $$i_1\circ(\tau^{A}_{X_1})^{-1}=(\tau^{A}_{P\otimes_EY_1\otimes_EP^\vee})^{-1}[1]\circ(A\otimes_Ai_1)\cma$$
and clearly
$$(A\otimes_Ai_1)\circ(g_2\otimes_AX_1)=(g_2\otimes_AP\otimes_EY_1\otimes_EP^\vee[1])\circ (P\otimes_EY_2\otimes_EP^\vee\otimes_A i_1)\fs$$
By the definition of $\psi$ we have
$$(P^\vee\otimes_Ai_1)\circ\psi=\eta^L_{Y_1\otimes_EP^\vee}[1]\circ(d_1\otimes_EP^\vee)$$
and so 
$$i_1\circ\bar{g}=(\tau^{A}_{PY_1P^\vee})^{-1}[1]\circ g_2PY_1P^\vee[1]\circ PY_2\eta^L_{Y_1P^\vee}[1]\circ PY_2d_1P^\vee \circ P\tau^{E^\op,\sigma_1}_{Y_2} P^\vee[n_1] $$
where again we have omitted the tensor product signs.
We use the description
$$g_2=\ev\circ(\tau^{E^\op}_{P})^{-1}\otimes_EP^\vee\circ P\otimes_Ef_2\otimes_EP^\vee$$
from Lemma \ref{altg2}.  Then
$$g_2PY_1P^\vee\circ PY_2\eta^L_{Y_1P^\vee}=\ev PY_1P^\vee\circ P\eta^L_{Y_1P^\vee}\circ(\tau^{E^\op}_P)^{-1}Y_1P^\vee\circ Pf_2Y_1P^\vee$$
by the naturality of $\eta^L$, and as $\ev PY_1 P^\vee=\tau^A_{PY_1P^\vee}\circ\varepsilon^L_{PY_1P^\vee}$ and $\varepsilon^L_{PY_1P^\vee}\circ P\eta^L_{Y_1P}=\id_{PY_1P^\vee}$ we can simplify this equality to obtain
$$g_2PY_1P^\vee\circ PY_2\eta^L_{Y_1P^\vee}=\tau^A_{PY_1P^\vee}\circ(\tau^{E^\op}_P)^{-1}Y_1P^\vee\circ Pf_2Y_1P^\vee\cma$$
so
\begin{align*}
i_1\circ\bar{g}&=(\tau^{A}_{PY_1P^\vee})^{-1}[1]
\circ \tau^A_{PY_1P^\vee}[1]\circ(\tau^{E^\op}_P)^{-1}Y_1P^\vee[1]\circ Pf_2Y_1P^\vee[1]
\circ PY_2d_1P^\vee \circ P\tau^{E^\op,\sigma_1}_{Y_2} P^\vee[n_1]\\
 &=(\tau^{E^\op}_P)^{-1}Y_1P^\vee[1]\circ Pf_2Y_1P^\vee[1]
\circ PY_2d_1P^\vee \circ P\tau^{E^\op,\sigma_1}_{Y_2} P^\vee[n_1]\fs
\end{align*}
we can rewrite this as
$$i_1\circ\bar{g}=(\tau^{E^\op}_P)^{-1}Y_1P^\vee[1]\circ PEd_1P^\vee\circ P(f_2E_{\sigma_1}\circ\tau^{E^\op,\sigma_1}_{Y_2}) P^\vee[n_1]\fs$$
Then by the naturality of $\tau$ and a simple commutation relation we have
$$i_1\circ\bar{g}=Pd_1P^\vee\circ  (\tau^{E^\op}_P)^{-1}E_{\sigma_1}P^\vee[n_1]\circ
 P\tau^{E^\op,\sigma_1}_{E}P^\vee[n_1] \circ P(f_2)_{\sigma_1}P^\vee[n_1]\fs$$
Finally, noting that $(\tau^{E^\op}_P)^{-1}E_{\sigma_1}P^\vee\circ P\tau^{E^\op,\sigma_1}_{E}P^\vee=\id_{PE_{\sigma_1}P^\vee}$, we get
$$i_1\circ\bar{g}=Pd_1P^\vee\circ P(f_2)_{\sigma_1}P^\vee[n_1]=-P\delta P^\vee$$
as required.

\textit{Step 5:}
By applying the octahedral axiom as outlined in Step 3, we obtain a distinguished triangle
$$X_2\otimes_AX_1\arr{}P\otimes_EY\otimes_EP^\vee[1]\arr{}A[1]\rsar{}$$
which we then rotate twice to get the triangle
$$P\otimes_EY\otimes_EP^\vee\arr{\hat{g}}A\arr{}X_2\otimes_AX_1\rsar{}\fs$$
The octahedral axiom also tells us that
$\hat{g}\circ P\otimes_E\iota\otimes_EP^\vee=-g_1\fs$

We wish to compare this triangle to the triangle $\nabla$ obtained in Step 2:
$$P\otimes_E Y\otimes_E P^\vee\arr{g} A\arr{h} X\rsar{i}$$
From the octahedron in Step 1 we see that $f\circ\iota=f_1$, so from the definition of $g$ in Step 2, $g\circ  P\otimes_E\iota\otimes_EP^\vee=g_1$.  Hence
$$(g+\hat{g}) P\otimes_E\iota\otimes_EP^\vee=0\in\Hom_{\Db(E\da E)}(P\otimes_E Y_2\otimes_E P^\vee,A)\fs$$
As $Y$ is a truncated resolution of $E$, Lemma \ref{mostproj} tells us that $P\otimes_EY\otimes_EP^\vee$ is a bounded complex of projective $A\da A$-bimodules, so $$\Hom_{\Db(A\da A)}(P\otimes_E Y_2\otimes_E P^\vee,A)=\Hom_{\Kb(A\da A)}(P\otimes_E Y_2\otimes_E P^\vee,A)$$
and we can work directly with complexes of $A\da A$-bimodules up to homotopy.  Then as $g$ and $\hat{g}$ are both nonzero only in degree zero, and $\iota$ is an isomorphism in degree zero, we have that $g+\hat{g}=0$ in the homotopy category and therefore $g=-\hat{g}$ in $\Db(A\da A)$.  So by the 5-lemma for triangulated categories we have an isomorphism of triangles
$$\xymatrix{
P\otimes_EY\otimes_EP^\vee\ar[r]^{\phantom{spacing}\hat{g}}\ar[d]^{-1} & A\ar[r]\ar@{=}[d] & X_2\otimes_AX_1\ar@{~>}[r]\ar@{-->}[d]^{\sim} & \\
P\otimes_E Y\otimes_E P^\vee\ar[r]^{\phantom{spacing}g} & A\ar[r]^{h} & X\ar@{~>}[r]^{i} &
}$$
which tells us that $X_2\otimes_AX_1\cong X$ as complexes of $A\da A$-bimodules.  Hence $\Psi_{P,Y}=X\otimes_A-$ and $\Psi_{P,Y_2}\circ \Psi_{P,Y_1}=X_2\otimes_AX_1\otimes_A-$ are naturally isomorphic.
\end{pf}
\end{prop}

\subsection{Inverses}
Recall that we are assuming that $P$ has no two isomorphic direct summands, so $P$ is a summand of $A$.  This will simplify the notation and, by Lemma \ref{nomult}, is no real restriction.

We consider quasi-inverses of the functors $\Psi_P$.  Given a truncated resolution $Y$ of $E$, we can apply to the distinguished triangle $\Delta$
$$Y\arr{f}E\arr{e} {E_\sigma}[n]\rsar{d}$$
the functor $-_{\sigma^{-1}}[-n]$ to obtain the triangle
$$Y_{\sigma^{-1}}[-n]\arr{f_{\sigma^{-1}}[-n]}E_{\sigma^{-1}}[-n]\arr{e_{\sigma^{-1}}[-n]} {E}\rsar{d_{\sigma^{-1}}[-n]}\fs$$
We write $Y'=Y_{\sigma^{-1}}[1-n]$, $f'=d_{\sigma^{-1}}[-n]$, $e'=e_{\sigma^{-1}}[-n]$, and $d'=f_{\sigma^{-1}}[-n]$ and so we have a distinguished triangle
$$Y'[-1]\arr{d'}E_{\sigma^{-1}}[-n]\arr{e'}E\rsar{f}$$
which we denote $\Delta'$.

As $P^\vee\otimes_A-$ is an exact functor, it has both a left and right adjoint.  By the tensor-hom adjunctions, we see that the right adjoint is $\RHom_E(P^\vee,-)$.  We obtain a unit $\eta':\id_A\to \RHom_E(P^\vee,P^\vee\otimes_A-)$ and a counit $\varepsilon':P^\vee\otimes_A\RHom_E(P^\vee,-)\to\id_E$.
We let
$$g'=\RHom_E(P^\vee,f'\otimes_EP^\vee)\circ \RHom_E(P^\vee,(\psi^L)^{-1})\circ \eta'_A:A\to \RHom_E(P^\vee,Y'\otimes_EP^\vee)$$
where $\psi^L$ is the isomorphism $\tau^{A^\op}_{P^\vee}\circ(\tau^E_{P^\vee})^{-1}:E\otimes_EP^\vee\arr{\sim}P^\vee\otimes_AA$.  Then we let
$X'=\cone(g)[-1]$ and so get a new distinguished triangle
$$X'\arr{}A\arr{g'} \RHom_E(P^\vee,Y'\otimes_EP^\vee)\rsar{}$$
which we denote $\nabla'$.  By arguments similar to those of Section \ref{pfofeq} we see that $X'\otimes_A-$ is an autoequivalence of $\Db(A)$ which sends $P$ to $P[-n]$ and acts as the identity on elements of $P^\perp$.  We denote this functor by $\Psi_{P,Y}'$, or simply $\Psi_P'$, and will show that it is a quasi-inverse of $\Psi_P$.

\begin{rmk}
As in Lemma \ref{altg2}, we could define $g$ using functors $-\otimes_AP$ and $\RHom_{E^\op}(P,-)$ and this would give a second unit-counit pair.  But we will not need to consider this alternative definition, and so write $\varepsilon'$ instead of $\varepsilon'^R$ to avoid cumbersome superscripts.
\end{rmk}

As $A$ is symmetric, it follows from Lemma \ref{perdual} that $E$ is also symmetric: we obtain an isomorphism of vector spaces $E\arr{\sim}E^*$ which is in fact a bimodule isomorphism by the naturality of Lemma \ref{perdual}.  Given a choice of isomophism
$$\varphi_A:A\arr{\sim}A^*$$
of $A\da A$-bimodules, we can construct an isomorphism
$$\varphi_E:E\arr{\sim}E^*$$
of $E\da E$-bimodules by the following rule.  Suppose $P=Aa_P$, so $a_P$ is the idempotent corresponding to $P$.  Then for $e,f\in E$, $\left(\varphi_E(e)\right)(f)$ is the image of $a_P$ in the following composition of maps:
$$P\arr{f}P\arr{e}P\arr{\iota_P}A\arr{\varphi_A(1_A)}k$$
where $\iota_P$ is the obvious injection. 
In particular, $\varphi_E(1_E)(1_E)=\varphi_A(1_A)(a_P)$. 

\begin{lem}\label{dualresn}
There is an isomorphism of triangles $\Delta'\cong \Delta^*$ such that the following diagram commutes:
$$\xymatrix{
E_{\sigma^{-1}}[-n]\ar[r]^{\phantom{spac}e'}\ar@{-->}[d] & E\ar[r]^{\phantom{s}f'}\ar[d]^{\varphi_E} & Y'\ar@{~>}[r]^(0.6){-d'[1]\phantom{sp}}\ar[d]^{i[-1]} & \\
(E_\sigma)^*[-n]\ar[r]^{\phantom{spac}e^*} & E^*\ar[r]^{\phantom{s}f^*} & Y^*\ar@{~>}[r]^(0.6){-d^*[1]\phantom{sp}} & 
}$$

\begin{pf}
As usual, we only need to prove that the right-hand square commutes, and the isomorphism of triangles will follow.  As $E$ is twisted-periodic it is self-injective, so $Y'$ and $Y^*$ are truncated injective resolutions of $E$ with kernel maps $f':E\into Y'$ and $f^*\circ \varphi_E:E\into Y^*$.  By duality, Proposition \ref{truncresisom} is also true when working with injective instead of projective resolutions, and by examining Lemma 2.1 of \cite{mannantrunc} we see that the square in our diagram commutes.
\end{pf}
\end{lem}

\begin{prop}
$\Psi_P'$ is a quasi-inverse of $\Psi_P$.

\begin{pf}
As $A$ is symmetric, we know that a quasi-inverse of $\Psi_P=X\otimes_A-$ is given by $X^*\otimes_A-$, and so our proof proceeds by showing that $X'$ and $X^*$ are isomorphic objects of $\Db(A\da A)$.

Consider the distinguished triangles $\nabla'$
$$X'\arr{}A\arr{g'} \RHom_E(P^\vee,Y'\otimes_EP^\vee)\rsar{}$$
and $\nabla^*$
$$X^*\arr{}A^*\arr{g^*}(P\otimes_EY\otimes_EP^\vee)^*\rsar{}\fs$$
It suffices to show that there exists an endomorphism $\hat{\alpha}$ of $P^\vee$ such that $\RHom_E(P^\vee,\hat{\alpha})$ is an isomorphism and the following diagram commutes:
$$\xymatrix{
A\ar[r]^{g'\phantom{spacing}}\ar[dddd]^{\varphi_A} & \RHom_E(P^\vee,Y'\otimes_EP^\vee)\ar[d]^{\RHom_E(P^\vee,i\otimes_E\hat{\alpha})} \\
 & \RHom_E(P^\vee,Y^*\otimes_EP^\vee)\ar[d]^{\RHom_E(P^\vee,Y^*\otimes\upsilon_P)} \\
 & \RHom_E(P^\vee,Y^*\otimes_EP^*)\ar[d]^{\RHom_E(P^\vee,\omega_{Y,P})} \\
 & \RHom_E(P^\vee,(P\otimes_EY)^*)\ar[d] \\
A^*\ar[r]^{g^*\phantom{spacing}} & (P\otimes_EY\otimes_EP^\vee)^*
}$$
Here, the isomorphism $i:Y'\to Y^*$ comes from Lemma \ref{dualresn}, $\omega_{M,N}^A$ is the natural isomorphism $M^*\otimes_AN^*\arr{\sim}(N\otimes_AM)^*$, $\upsilon_M$ is the natural isomorphism $M^\vee\arr{\sim} M^*$, and the bottom vertical map is given by the tensor-hom adjunction between $P\otimes_EY\otimes_E-$ and $\Hom_k(P\otimes_EY,-)$.  All vertical maps in this diagram are isomorphisms, so after showing that the diagram commutes the statement of the proposition will follow from the 5-lemma for triangulated categories.

Recall that
$$g=\varepsilon^L_A\circ(P\otimes_E\psi^L)\circ(P\otimes_Ef\otimes_EP^\vee)$$
so
$$g^*=(P\otimes_Ef\otimes_EP^\vee)^*\circ(P\otimes_E\psi^L)^*\circ(\varepsilon^L_A)^*\cma$$
and
$$g'=\RHom_E(P^\vee,f'\otimes_EP^\vee)\circ \RHom_E(P^\vee,(\psi^L)^{-1})\circ \eta'_A$$
where $f'=d_{\sigma^{-1}}[-n]$.  

Let $\psi^R$ be the isomorphism $\tau^A_P\circ(\tau^{E^\op}_P)^{-1}:P\otimes_EE\arr{\sim}A\otimes_AP$.  Then we have the following diagram, where the composition of the vertical maps on the right hand side is just the composition of the maps in the top right hand sides of the above square:
$$\xymatrix{
 A^*\ar[d]^{\eta'_{A^*}} & A\ar[d]_{\eta'_A}\ar[l]_{\varphi_A} \\
 \RHom_E(P^\vee,P^\vee\otimes_AA^*)\ar[d]_{\RHom_E(P^\vee,\omega^A_{P,A}\circ(\upsilon_P\otimes A^*))} & \RHom_E(P^\vee,P^\vee\otimes_AA)\ar[d]^{\RHom_E(P^\vee,(\psi^L)^{-1})}\ar[l]_{P^\vee\otimes\varphi_A} \\ 
 \RHom_E(P^\vee,(A\otimes_AP)^*)\ar[d]_{\RHom_E(P^\vee,(\psi^R)^*)} & \RHom_E(P^\vee,E\otimes_EP^\vee)\ar[d]^{\RHom_E(P^\vee,f'\otimes_EP^\vee)}\ar[ldd]_(0.4){\varphi_E\otimes \hat{\alpha}} \\
 \RHom_E(P^\vee,(P\otimes_EE)^*)\ar[d]_{\RHom_E(P^\vee,(E^*\otimes\upsilon_P^{-1})\circ\omega_{E,P}^{-1})}& \RHom_E(P^\vee,Y'\otimes_EP^\vee)\ar[d]^{\RHom_E(P^\vee,i\otimes_E\hat{\alpha})} \\
 \RHom_E(P^\vee,E^*\otimes_EP^\vee)\ar[r]^{f^*\otimes P^\vee}\ar[d]_{\RHom_E(P^\vee,\omega^E_{E,P}\circ(E^* \otimes\upsilon_P))} & \RHom_E(P^\vee,Y^*\otimes_EP^\vee)\ar[d]^{\RHom_E(P^\vee,\omega^E_{Y,P}\circ(Y^* \otimes\upsilon_P))} \\
 \RHom_E(P^\vee,(P\otimes_EE)^*)\ar[r]^{(P\otimes f)^*}\ar[d] & \RHom_E(P^\vee,(P\otimes_EY)^*)\ar[d] \\
 (P\otimes_EE\otimes_EP^\vee)^*\ar[r]^{(P\otimes f\otimes P^\vee)^*} & (P\otimes_EY\otimes_EP^\vee)^*
}$$
We have omitted the functor $\RHom_E(P^\vee,-)$ from the middle horizontal maps to save space.

The four squares in the diagram commute because of the naturality of $\eta'$, the relation $i f'=f^* \varphi_E$ from Lemma \ref{dualresn}, the naturality of $\omega$, and the naturality of the tensor-hom adjunctions.

To show that the hexagon in the diagram commutes, we consider the composition
$$P^\vee\arr{}P^*\to P^*\otimes_AA\arr{} P^*\otimes_AA^*\arr{}(A\otimes_AP)^*\arr{}P^*\arr{} P^\vee$$
which just acts as multiplication by the constant $\lambda_A=\varphi_A(1_A)(1_A)\in k$.  A similar statement holds for tensoring on the left by $E$, with a constant $\lambda_E$.  So by the commutativity of the diagram of isomorphisms
$$\xymatrix@C=2pt{
P^*\otimes_AA^*\ar[rr]^{\omega^A_{P,A}} & & (A\otimes_AP)^*\ar[rr]^{(\psi^R)^*}\ar[rd]_(0.3){(\tau^A_P)^*} & & (P\otimes_EE)^*\ar[ld]^(0.3){(\tau^{E^\op}_P)^*} & & E^*\otimes_EP^*\ar[ll]_{\omega^E_{E,P}}\\
 & P^*\otimes_AA\ar[ul]^(0.6){P^*\otimes\varphi_A} && P^*\ar[ll]^(0.4){\tau^{A^\op}/\lambda_A}\ar[rr]_(0.4){\tau^E_{P^*}/\lambda_E} && E\otimes_EP^*\ar[ur]_(0.6){\varphi_E\otimes P^*}
}$$
the hexagon in our large diagram commutes after an appropriate rescaling by $\lambda=\lambda_A/\lambda_E$.

We note that by naturality, the following diagram also commutes:
$$\xymatrix{
\RHom_E(P^\vee,(A\otimes_AP)^*)\ar[d]_{\RHom_E(P^\vee,(\psi^R)^*)}\ar[r] & (A\otimes_AP\otimes_EP^\vee)^*\ar[d]^{(\psi^R\otimes_EP^\vee)^*} \\
\RHom_E(P^\vee,(P\otimes_EE)^*)\ar[r] & (P\otimes_EE\otimes_EP^\vee)^*
}$$
and so to prove that the original square commutes, it only remains to show that the following diagram commutes:
$$\xymatrix{
 A^*\ar[d]^{(\varepsilon^R_A)^*}\ar[r]^(0.3){\eta'_{A^*}} & \RHom_E(P^\vee,P^\vee\otimes_AA^*)\ar[d]^{\RHom_E(\hat{\alpha},\omega_{P,A})} \\
 (A\otimes_AP\otimes_EP^\vee)^* & \RHom_E(P^\vee,(A\otimes_AP)^*)\ar[l]
}$$
To see this, we will use the fact \cite[Theorem IV.1]{mac} that $\eta'_{A^*}$ is a universal arrow from $A^*$ to $\RHom_E(P^\vee,-)$, and that $(\varepsilon_A^R)^*$ is a universal arrow from $A^*$ to $(-\otimes_EP^\vee)^*$.  We denote by $\gamma$ the composition of the isomorphisms
$$\RHom_E(P^\vee,P^\vee\otimes_AA^*)\arr{\sim}
\RHom_E(P^\vee,(A\otimes_AP)^*)\arr{\sim}(A\otimes_AP\otimes_EP^\vee)^*\fs$$
Because $\gamma^{-1}\circ (\varepsilon^R_A)^*$ is a map from $A^*$ to $\RHom_E(P^\vee,P^\vee\otimes_AA^*)$ there is some endomorphism $\alpha$ of $P^\vee\otimes_AA^*$ such that $\gamma^{-1}\circ (\varepsilon^R_A)^*=\RHom_E(P^\vee,\alpha)\circ\eta'_{A^*}$, and similarly there is some endomorphism $\beta$ of $A\otimes_AP$ such that $\gamma\circ\eta'_{A^*}=(\beta\otimes_EP^\vee)^*\circ(\varepsilon^R_A)^*$.
Combining the two equations, we obtain
$$(\varepsilon^R_A)^*=\gamma\circ\RHom_E(P^\vee,\alpha)\circ\gamma^{-1}\circ(\beta\otimes_EP^\vee)^*\circ(\varepsilon^R_A)^*$$
and so by the universal property of $(\varepsilon^R_A)^*$ we see that
$$\gamma\circ\RHom_E(P^\vee,\alpha)\circ\gamma^{-1}\circ(\beta\otimes_EP^\vee)^*=\id\fs$$
Hence $\RHom_E(P^\vee,\alpha)$ is an epimorphism and $(\beta\otimes_EP^\vee)^*$ is a monomorphism and, by the universal property of $\eta'_{A^*}$, vice versa.  We define $\hat{\alpha}$ as the composition
$$P^\vee\arr{\sim}P^\vee\otimes_AA^*\arr{\alpha}P^\vee\otimes_AA^*\arr{\sim}P^\vee$$
so $\RHom_E(P^\vee,Y'\otimes_E\hat\alpha)$ is an isomorphism by naturality.  Then, by construction, all our diagrams commute.
\end{pf}
\end{prop}

\section{Deconstructing autoequivalences}

We will show that the derived autoequivalences described above are actually compositions of derived equivalences which are already known.  This will produce a ``circle of equivalences'': a chain of derived equivalent algebras starting and ending with $A$.

\subsection{Tilting theory}
We recall the main result on equivalences of derived categories, due to Rickard.

\begin{defn}[\cite{ricmor}]
A bounded complex $T$ of projective $A$-modules is a \emph{one-sided tilting complex} (or just a \emph{tilting complex}) if
\begin{enumerate}
\item $\Hom_{\Kb(A)}(T,T[m])=0$ for all $m\neq0$, and
\item $T\add$ generates $\Kb(A)$, i.e., $\Kb(A)$ has no proper full triangulated subcategory, closed under isomorphisms and direct summands, which contains $T\add$.
\end{enumerate}
\end{defn}

\begin{thm}[\cite{ricmor}]\label{dermor}
For two algebras $A$ and $B$, $\Db(A)$ and $\Db(B)$ are equivalent as triangulated categories if and only if $B\cong\End_{\Kb(A)}(T)$ where $T$ is a tilting complex.
If such a tilting complex exists, there is an equivalence $\Db(A)\to \Db(B)$ which sends the tilting complex $T$ to the regular $B$-module. 
\end{thm}

The previous theorem gives a good way to produce derived equivalences.  However, it can be difficult to calculate the image of an arbitrary complex under a functor which exists by the theorem.  Calculations are made easier by the use of \emph{two-sided tilting complexes}.  These are complexes $X\in\Db(A\da B)$ such that $\RHom_A(X,-):\Db(A)\to\Db(B)$, or equivalently $X\dert_B-:\Db(B)\to\Db(A)$, is an equivalence of triangulated categories.

\begin{prop}[\cite{rictwo}]\label{exists2sided}
If $F:\Db(A)\arr{\sim}\Db(B)$ is a derived equivalence then there exists a two-sided tilting complex $X\in\Db(B\da A)$ such that $FV\cong X\dert_AV$ for all $V\in\Db(A)$ and $FP=X\otimes_AP$ for all projective $A$-modules $P$.
\end{prop}

\subsection{Combinatorial tilting complexes}\label{combsec}
We present the construction of some well-known tilting complexes and some of their basic properties. Most, if not all, of the material in this subsection is well known, but in some cases it can be difficult to obtain suitable references so we include details for the convenience of the reader.

As always, $A$ is a symmetric algebra.  Let $\{P_i\}_{i=1}^r$ be representatives of a complete set of isomorphism classes of indecomposable projective $A$-modules and $\{S_i\}_{i=1}^r$ be a corresponding set of simple modules.  Let $J$ be a nontrivial proper subset of $I=\{1,\ldots,r\}$ and write $P=\bigoplus_{j\in J}P_j$.
\begin{defn}\label{combtilt}
A \emph{combinatorial tilting complex at $J$} is a complex $T=\bigoplus_{i\in I}T_i^{\ell_i}$, with $\ell_i\geq1$ for all $1\leq i\leq r$, concentrated in degrees $-1$ and $0$ where
$$T_i = \left\{ \begin{array}{ll}
P_i\to 0 & \textrm{if }i\in J\\
P'_i\arr{\varphi_i} P_i & \textrm{if }i\in I\bs J
\end{array} \right.
$$
such that for each $i\in I\bs J$,
\begin{enumerate}
\item $P'_i\in P\add$,
\item $\coker\varphi_i$ has no composition factors isomorphic to $S_j$ for $j\in J$, and
\item $P'_i$ is minimal with respect to its direct sum decomposition.
\end{enumerate}
%
%
If $\ell_i=1$ for all $1\leq i\leq r$ then we say $T$ is \emph{basic}.
\end{defn}
Combinatorial tilting complexes were first introduced by Rickard  in the case where $J$ consists of a single element \cite[Example 3.III.4]{ricthesis}, and Okuyama noted that the construction could be immediately generalised to arbitrary subsets $J$ of $I$ \cite{okuy}.  They have appeared under many names, such as ``special $2$-term tilting complexes'' and ``Okuyama-Rickard tilting complexes'', and are well studied: see for example \cite{okuy,kztilt,jsz,ai}.

To show that combinatorial tilting complexes exist, we give the following explicit construction.
\begin{lem}
For each $i\in I\bs J$, let $Q$ be the intersection of all submodules $Q'$ of $P_i$ such that the inclusion $Q'\into P_i$ induces an isomorphism
$$\Hom_A(P,Q')\arr{\sim}\Hom_A(P,P_i)$$
and let $\pi_Q:P_Q\to Q$ be a projective cover of $Q$.  Then if we let $P'_i=P_Q$ and let $\varphi_i$ be the composition
$$P_Q\onto Q\into P_i$$
the modules $P_i'$ and maps $\varphi_i$ define a combinatorial tilting complex at $J$.

\begin{pf}
First we show that $P'_i$ belongs to $P\add$.  Suppose for contradiction that it has a summand $P_m$ for $m\in I\bs J$.  Then the head of $Q$ must contain a summand isomorphic to $S_m$, and the kernel $K$ of the projection $Q\onto S_m$ is a proper submodule of $Q$ such that the inclusion $K\into P_i$ induces an isomorphism $\Hom_A(P,K)\arr{\sim}\Hom_A(P,P_i)$,
contradicting the definition of $Q$.

To see that $\coker\varphi_i$ has no composition factors isomorphic to $S_j$ for $j\in J$ we apply the exact functor $\Hom_A(P,-)$ to the exact sequence
$$0\to Q\into P_i\onto\coker\varphi_i\to 0$$
and as the map $Q\into P_i$ is sent to an isomorphism, $\Hom_A(P,\coker\varphi_i)=0$.

To show that $P'_i$ is minimal with respect to direct sum decomposition, we suppose that there exists some $P''\in P\add$ and a map $\psi:P''\to P_i$ such that $\coker \psi$ has no composition factors isomorphic to $S_j$ for $j\in J$, and we will show that $P'_i$ is isomorphic to a summand of $P''$.  First, we apply $\Hom_A(P,-)$ to the short exact sequence
$$0\to\im\psi\into P_i\onto\coker\psi\to0$$
to deduce that $\im\psi\into P_i$ induces an isomorphism $$\Hom_A(P,\im\psi)\arr{\sim}\Hom_A(P,P_i)\cma$$ and so $Q\subset\im\psi$.  Then as the composition
$$Q\into\im\psi\into P_i$$
becomes an isomorphism on application of $\Hom_A(P,-)$, the cokernel $L$ of the map $Q\into\im\psi$ has a composition series consisting entirely of simple modules $S_j$ for $j\in J$.  But then as $L$ is a quotient of $P''\in P\add$, $L$ must be the zero module, and so $Q\cong\im\psi$.  Therefore as $\im\psi$ is a quotient of $P''$, $P'_i=P_Q$ must be isomorphic to a summand of $P'_i$ by the minimality of projective covers.
\end{pf}
\end{lem}

We show that the name of these complexes is reasonable.
\begin{prop}
Combinatorial tilting complexes are tilting complexes.

\begin{pf}
Any combinatorial tilting complex $T$ is by definition a bounded complex of projective modules, so we just need to show that it generates $\Kb(A)$ and has no self-extensions.  As $T$ is concentrated in two adjacent degrees, by Lemma \ref{perdual} we only need to show that $\Hom_{\Kb(A)}(T,T[1])=0$.
Note that for any map $\varphi_i:P_i'\to P_i$, $P_i'$ must be in $P\add$ because, if $Q=\im\varphi_i$, by the minimality of $Q$, the semisimple module $Q/\rad Q$ cannot have a summand isomorphic to $S_i$ for $i\in I\bs J$, and the projective covers of $Q$ and $Q/\rad Q$ coincide.  So the only possible nonzero maps are of the form
$$\xymatrix{
 0\ar[r] & P'\ar[d]^f\ar[r] & \ldots \\
 P'_i\ar[r]^{\varphi_i} & P_i\ar[r] & 0
}$$
with $P'\in P\add$.  But as $\coker\varphi_i$ has no composition factors isomorphic to $S_j$ for $j\in J$, $f$ must factor through $\varphi_i$ and so this map is homotopic to zero.

To see that $T$ generates $\Kb(A)$, note that for all $j\in J$, $P_j\in T\add$ and so becuase we have distinguished triangles
$$P_i'\to P_i\to T_i\rsa$$
and $P_i'\in T\add$, $P_i$ must be in the triangulated category generated by $T$.  So all the isoclasses of indecomposable projectives are in this subcategory and hence it generates $\K^b(A)$.
\end{pf}
\end{prop}

\begin{cor}
For fixed $J\subset I$, all basic tilting complexes at $J$ are unique up to isomorphism.

\begin{pf}
As basic tilting complexes are tilting, this follows from Corollary 4.8 of \cite{jsz}, but we also give a direct proof.

Fix $i\in I\bs J$.  We only need to show that if the $A$-module $P'_i$ and the map $\varphi_i:P'_i\to P_i$ satisfy the conditions of Definition \ref{combtilt}, then this map is isomorphic to the map $P_Q\to P_i$ constructed in the previous lemma.  Clearly $\im\varphi_i$ is a submodule of $Q$.  Then given the following diagram
$$\xymatrix@=6pt{
P_i'\ar[rr]\ar@{->>}[dr]\ar@{-->}[dd] & &P_i\ar@{-->}[dd]\\
 & \im\varphi_i\ar@{^{(}->}[ur]\ar@{^{(}->}[dd] &\\
P_Q\ar'[r][rr]\ar@{->>}[dr] & &P_i\\
 & Q\ar@{^{(}->}[ur] &
}$$
the left dotted arrow exists by the projectivity of $P'_i$ and is an isomorphism by the minimality condition.  As $A$ is symmetric, $P_i$ is an injective module, so the right dotted arrow exists, and by commutativity it is an injective endomorphism of a finite-dimensional module and hence an isomorphism.
\end{pf}
\end{cor}

So by Theorem \ref{dermor} each combinatorial tilting complex $T$ induces a derived equivalence $$F_J:\Db(A)\arr{\sim}\Db(A^{(1)})\cma$$ where $A^{(1)}=\End_{\Db(A)}(T)^\op$, which we call the \emph{combinatorial tilt at $J$}.  We have a natural indexing of the projective $A^{(1)}$-modules $P^{(1)}_i=F_J(T_i)$ and the simple $A^{(1)}$-modules $S^{(1)}_i:=P^{(1)}_i/\rad P^{(1)}_i$ by $I$.

\begin{prop}
Combinatorial tilts are perverse, with filtrations of the simple modules of both $A$ and $B=\End_{\Db(A)}(T)$ given by the filtration
$$\emptyset\subset I\bs J\subset I$$
and the perversity function $p:\{1,2\}\to\Z$ given by $p(1)=0$ and $p(2)=-1$.

\begin{pf}
The proof is similar to that of Proposition \ref{pertwperv}, using the fact that $F_T(P)=P^{(1)}[-1]$, which follows from the definition.
\end{pf}
\end{prop}

Recall that, for a subcategory $\CC$ of a triangulated category $\T$, a map $f:V\to W$ is a \emph{right $\CC$-approximation of $W$} if $V\in C$ and $\Hom_\CC(C,f)$ is surjective for every $C\in\CC$.  The map $f:V\to W$ is \emph{right minimal} if all maps $g:V\to V$ such that $fg=f$ are isomorphisms, and is a \emph{right minimal approximation} if it is both right minimal and a right approximation.
\begin{lem}\label{combapprox}
The maps $\varphi_i$ in a combinatorial tilting complex at $J$ are right minimal $P\add$-approximations of $P_i$.

\begin{pf}
First we assume that $\varphi_i:P'_i\to P_i$ is a right minimal $P\add$-approximation of $P_i$.  Then $\Hom_A(P,\varphi_i)$ is surjective, and so has zero cokernel.  As $\Hom_A(P,-)$ is exact, $\coker\Hom_A(P,\varphi_i)=\Hom_A(P,\coker\varphi_i)$ and so $\coker\varphi_i$ has no composition factor isomorphic to $S_j$ for $j\in J$.  Now, to show minimality, let $Q=\varphi_i$ and let $R$ be a submodule of $Q$ such that $P_i/R$ has no composition factor isomorphic to $S_j$ for $j\in J$.  Then  $Q/R$ is a quotient of a module in $P\add$ such that $\Hom_A(P,Q/R)=0$ and hence $Q/R=0$.

To show that $\varphi_i$ is a projective cover of its image $Q$ we only need to show that $P_i'$ is minimal with respect to its diret sum decomposition.  But if $P_i'\cong M\oplus N$ where $Q$ is a quotient of $M$, then $\varphi_i$ is equal to a composition of morphisms $P_i'\onto M\into P_i'\onto Q$ and so by the minimality of the approximation $\varphi_i$ we have $P_i'\cong M$.

Suppose now that $\varphi_i:P'_i\to P_i$ satisfies the conditions of Definition \ref{combtilt}.  Reversing the argument above we can show that $\Hom(P,f)$ is surjective, and from this it follows that $f$ is a right $P\add$-approximation.  The right minimality of $f$ follows because if $g:P'_i\to P'_i$ is not an isomorphism it must eventually annihilate some direct summand of $P'_i$, contradicting the minimality condition of the definition.
\end{pf}
\end{lem}
In more modern language, this lemma tells us that, up to shift, combinatorial tilts at $J$ are just tilting mutations at $I\bs J$.

\subsection{Circles of equivalences}
Given a fixed subset $J$ of the indexing set $I$ of isoclasses of projective $A$-modules, the combinatorial tilt $F_J$ at $J$ gives us a new algebra which we denote $A^{(1)}$.  For the tilting complex $T$, the projective $A^{(1)}$-modules $P_i^{(1)}=F_J(T_i)$ are again indexed by $I$.  We can repeat this procedure to obtain another derived equivalence $\Db(A^{(1)})\arr{\sim}\Db(A^{(2)})$, which we also denote $F_J$, and continuing in this manner we define the \emph{$n$th iterated combinatorial tilt at $J$}, $$F_J^n:\Db(A^{(0)})\arr{\sim}\Db(A^{(n)})\cma$$ where $A^{(0)}=A$.  It is convenient to write $P^{(n)}=\bigoplus_{j\in J}P_j^{(n)}$.

It is easy to see that the $n$th iterated combinatorial tilt is perverse with filtration $\emptyset\subset I\bs J\subset I$
and perversity function $p:\{1,2\}\to\Z$ given by $p(1)=0$ and $p(2)=-n$.  It is natural to ask whether our periodic twists agree with the inverse of $n$th iterated combinatorial tilts, and we will see that they do.

We write $\Omega_A(M)$ for the kernel of a projective cover $P_M\onto M$ of the $A$-module $M$.
\begin{lem}\label{iterapprox}
For $i\in I\bs J$ and $n\in\Z_{>0}$, given a complex $T_i$ of the form
$$P'_{i,n}\arr{}P'_{i,n-1}\arr{}\ldots\arr{}P'_{i,1}\arr{}P_i$$
with $P'_{i,m}\in P\add$ for all $1\leq m\leq n$, $T_i$ is the preimage of $P_i^{(n)}$ under an $n$th iterated combinatorial tilt if and only if
$$\Hom_A(P,T_i)\cong \Omega^n_E\Hom_A(P,P_i)[n]$$
in $\Db(E)$.

\begin{pf}
We proceed by induction with Lemma \ref{combapprox} providing our base case.  First we suppose that $F_J^n(T_i)\cong P^{(n)}_i$.  Then we have $F_J^{n-1}(F_JT_i)\cong P^{(n)}$, so our inductive assumption tells us that 
$$\Hom_{A^{(1)}}(P^{(1)},F_JT_i)\cong \Omega_E^{n-1}\Hom_{A^{(1)}}(P^{(1)},P_i^{(1)})[n-1]\fs$$
As $P[n]$ and $T_i$ are both summands of the same tilting complex, we have an isomorphism of $E$-modules
$$\Hom_A(P,T_i)\cong \Hom_{\Db(A)}(P[n],T_i)[n]\cong \Hom_{\Db(A^{(1)})}(P^{(1)}[n-1],F_JT_i)[n]$$
and so we see that $\Hom_A(P,T_i)\cong \Omega_E^{n-1}\Hom_{A^{(1)}}(P^{(1)},P_i^{(1)})[n]$.  We use Lemma \ref{combapprox} to deduce that
\begin{align*}
\Hom_{A^{(1)}}(P^{(1)},P_i^{(1)}) &\cong \Hom_{\Db(A^{(1)})}(P^{(1)},P_i^{(1)})\\
  &\cong \Hom_{\Db(A)}(P[1],F_J^{-1}P_i^{(1)})\\
  &\cong \Hom_A(P,F_J^{-1}P_i^{(1)})[-1]\\
  &\cong \Omega_E\Hom_A(P,P_i)
\end{align*}
and so we have $\Hom_A(P,T_i)\cong \Omega_E^{n}\Hom_A(P,P_i)[n]$, as required.

Now suppose that $T_i$ has the form
$$P'_{i,n}\arr{}P'_{i,n-1}\arr{}\ldots\arr{}P'_{i,1}\arr{}P_i$$
and that $\Hom_A(P,T_i)\cong \Omega^n_E\Hom_A(P,P_i)[n]$.  If we can show that $F_JT_i$ has the form
$$P_{i,n-1}''\arr{}P''_{i,n-2}\arr{}\ldots\arr{}P''_{i,1}\arr{}P^{(1)}_i$$
with $P''_{i,m}\in P^{(1)}\add$, 
and that $\Hom_{A^{(1)}}(P^{(1)},F_JT_i)\cong \Omega_E^{n-1}\Hom_{A^{(1)}}(P^{(1)},P_i^{(1)})$, then our inductive assumption will tell us that $F^{n-1}_J(F_JT_i)\cong P_i^{(n)}$, which is what we want.

To see that $F_JT_i$ has the desired form, we note that as $\Hom_A(P,T_i)$ is concentrated in degree $n$, the map obtained by applying $\Hom_A(P,-)$ to $P'_{i,1}\arr{}P_i$ must be surjective with $P''_{i,1}$ minimal with respect to direct sum decomposition.  So $F_J$ sends $P'_{i,1}\arr{}P_i$ to $P_i^{(1)}$ and $P'_{i,m}\in P\add$ to $P'_{i,m-1}[-1]$, where $P'_{i,m-1}\in P^{(1)}\add$.

Finally we show that $\Hom_{A^{(1)}}(P^{(1)},F_JT_i)\cong \Omega_E^{n-1}\Hom_{A^{(1)}}(P^{(1)},P_i^{(1)})$.  Using our assumption that $\Hom_A(P,T_i)\cong \Omega^n_E\Hom_A(P,P_i)[n]$, we see
$$\Hom_{A^{(1)}}(P^{(1)},F_JT_i) \cong \Hom_A(P[1],T_i)
  \cong \Omega^n_E\Hom_A(P,P_i)[n-1]$$
and, as above, $\Hom_{A^{(1)}}(P^{(1)},P_i^{(1)})\cong \Omega_E\Hom_A(P,P_i)$, so we are done.
\end{pf}
\end{lem}

\begin{thm}\label{circles}
If $P$ is a projective module for a symmetric algebra $A$ and $E=\End_A(P)^\op$ is twisted periodic of period $n$, then the periodic twist $\Psi_P$ agrees with the inverse $F^{-n}_J$ of an $n$th iterated combinatorial twist, i.e., for all $V\in \Db(A)$, $X\otimes_AV\cong F_J^{-n}V$.  In particular, the algebra $A^{(n)}$ is isomorphic to $A$.

\begin{pf}
Without loss of generality we will assume that $A$ is a basic algebra: if this is not the case, the proof can be adapted by increasing the multiplicity of summands of our combinatorial tilting complexes.  By Lemma \ref{nomult}, we may also assume that $P=\bigoplus_{j\in J}P_j$.

First we will show that $\Psi_P$ and $F^{-n}_J$ agree on the regular $A$-module $A=\bigoplus_{i\in I}P_i$. 
By Corollary \ref{pn} and the definition of combinatorial twists, both functors send $P$ to $P[n]$.  For projective modules $P_i$ with $i\in I\bs J$, consider $\Psi_P(P_i)=X\otimes_AP_i$.  Using the distinguished triangle $\nabla$ we see that $X\otimes_AP_i$ has the form described in Lemma \ref{iterapprox}, so if we can show that $\Hom_A(P,X\otimes_AP_i)\cong \Omega^n_E\Hom_A(P,P_i)[n]$ the functors must agree on $P_i$.  

As $E$ is $\sigma$-periodic of period $n$ for some automorphism $\sigma$ of $A$, $$\Omega^n_E\Hom_A(P,P_i)\cong E_\sigma\otimes_E\Hom_A(P,P_i)\fs$$  By the isomorphism of triangles $\Delta\otimes_EP^\vee\cong P^\vee\otimes_A\nabla$, proved similarly to Proposition \ref{upispd}, we see that $P^\vee\otimes_AX\cong E_\sigma\otimes_EP^\vee[n]$.  Then $$\Hom_A(P,X\otimes_AP_i)\cong P^\vee\otimes_AX\otimes_AP_i\cong E_\sigma\otimes_EP^\vee\otimes P_i[n]\cong E_\sigma\otimes_E\Hom_A(P,P_i)\cma$$ so $\Hom_A(P,X\otimes_AP_i)\cong \Omega^n_E\Hom_A(P,P_i)[n]$ as required.

To see that both functors agree on arbitrary $V\in\Db(A)$, use Proposition \ref{exists2sided} to obtain a two-sided tilting complex $W$ that agrees with $F^{-n}_J$ on all objects of $\Db(A)$, and then as
$$W\otimes_AX^*\otimes_AA\cong A \cong X^*\otimes_AW\otimes_AA$$
we must have $W\cong X$ in $\Db(A\da A)$.

Finally, as $F_J^{-n}(A)\cong X$, we have isomorphisms of algebras
$$A^{(n)}=\End_{\Db(A)}(X)^\op\cong\Hom_{\Db(A)}(A,X^*\otimes_AX\otimes_AA)^\op\cong A\fs$$
\end{pf}
\end{thm}

This chain of equivalences can be pictured as follows:
$$\xymatrix  @R=24pt @C=10pt  {
   & \Db(A)\ar[dr]^{F_J} & \\
 \Db(A^{(n-1)})\ar[ur]^{F_J} & & \Db(A^{(1)})\ar[d]^{F_J}\\
 \Db(A^{(n-2)})\ar[u]^{F_J} & & \Db(A^{(2)})\ar[dl] \\
 & {}\phantom{n}\:\hdots\:\phantom{n}\ar[ul] & 
}$$

As each $F_J$ is a derived equivalence, travelling $n$ steps around the circle will always give an autoequivalence.  By considering the image of our projective $A$-module $P$ in these combinatorial tilts, we can describe each such autoequivalence.
\begin{lem}
Given a diagram as above and any $0\leq m\leq n-1$, the composition $$(F_J)^{-n}:\Db(A^{(m)})\to\Db(A^{(m)})$$ agrees with the periodic twist $\Psi_{P^{(m)}}$.
\end{lem}

\section{Examples}
\subsection{Truncated polynomial rings}\label{truncex}
Consider the commutative algebra $E=k[x]/\gen{x^{n+1}}$ for $n\geq1$.  $E$ is periodic, as can be seen from the following exact sequence of $E\da E$-bimodules:
$$\ldots\to E\otimes_kE\arr{d_3}E\otimes_kE\arr{d_2}E\otimes_kE\arr{d_1}E\otimes_kE\arr{m} E\to 0$$
where, for $i$ even, $d_i=\sum_{j=0}^n x^j\otimes x^{n-j}$, for $i$ odd, $d_i=1\otimes x-x\otimes 1$, and $m$ is the multiplication map $e\otimes e'\mapsto ee'$.  This is a special case of Proposition 1.3 of \cite{bach}.  The resolution is periodic with periodicity $2$, and so we have a short exact sequence
$$0\to E[1]\into Y\onto E\to0$$
of chain complexes of $E\da E$-bimodules where $Y$ is the complex
$$E\otimes_kE\stackrel{1\otimes x-x\otimes 1}{\longrightarrow}E\otimes_kE$$
concentrated in degrees $-1$ and $0$.

Suppose that $A$ is a symmetric algebra with a projective module $P$ such that $E=\End_A(P)\cong k[x]/\gen{x^{n+1}}$.  Using the above data, we can construct the periodic twist $\Psi_P=X\otimes_A-$.  Using the isomorphism $P^*\cong P^\vee$, we see that the complex of $A\da A$-bimodules $X$ is isomorphic to the following complex
$$P\otimes_kP^\vee\stackrel{1\otimes x-x\otimes 1}{\longrightarrow}P\otimes_kP^\vee\stackrel{\ev}{\longrightarrow}A$$
where $\ev$ is the obvious evaluation map $P\otimes_k\Hom_A(P,A)\to A$ sending $p\otimes f$ to $f(p)$.  Using the isomorphism of functors $P^\vee\otimes_A-\cong \Hom_A(P,-)$ we see that we have constructed a complex of functors
$$P\otimes_k\Hom_A(P,-)\stackrel{1\otimes x-x\otimes 1}{\longrightarrow}P\otimes_k\Hom_A(P,-)\stackrel{\ev}{\longrightarrow}-$$
which is an analogue for symmetric algebras of the $\BP^{n}$-twists defined by Huybrechts and Thomas \cite{ht}.

We now specialise to the case $n=1$.  Then $E=k[x]/\gen{x^2}$ has the following resolution in terms of $E\da E$-bimodules:
$$\ldots\longrightarrow E\otimes_kE\stackrel{1\otimes x-x\otimes 1}{\longrightarrow}E\otimes_kE\stackrel{1\otimes x+x\otimes 1}{\longrightarrow}E\otimes_kE\stackrel{1\otimes x-x\otimes 1}{\longrightarrow}E\otimes_kE$$
so in fact we have a short exact sequence of $E\da E$-bimodules
$$0\to E_\sigma\stackrel{\iota}{\into} E\otimes_kE\stackrel{m}{\onto} E\to0$$
where $m$ is the multiplication map, $\iota$ sends $e$ to $e\otimes x-ex\otimes1$, and $\sigma$ is the automorphism of $E$ sending $x$ to $-x$.  Hence we see that $E$ is twisted-periodic of period $1$.  Similarly to the case above, if $A$ is a symmetric algebra with a projective module $P$ such that $E=\End_A(P)\cong k[x]/\gen{x^2}$, we obtain a two-sided tilting complex
$$P\otimes_kP^\vee\stackrel{\ev}{\longrightarrow}A$$
that gives a complex of functors
$$P\otimes_k\Hom_A(P,-)\stackrel{\ev}{\longrightarrow}-$$
and this gives an algebraic analogue of the spherical twists of Seidel and Thomas \cite{st}.  This case also recovers the autoequivalences constructed by Rouquier and Zimmermann \cite{rz} for Brauer tree algebras with no exceptional vertex, as such algebras satisfy $\End_A(P)\cong k[x]/\gen{x^2}$ for every indecomposable projective $P$.

Huybrechts and Thomas noted \cite{ht} that $\BP^1$-twists were just the square of spherical twists.  By direct calculation, the same phenomenon can be seen to hold with our algebraic analogues of these twists described above, but we can also provide a theoretical explanation.  By combining the period $1$ truncated resolutions of $E$ and $E_\sigma$, we obtain a period $2$ truncated resolution of $E$ as follows:
$$\xymatrix@=8pt{
0\ar[rr] &&E\ar[rr]^{\iota_\sigma} && Y_\sigma\ar[rd]_{m_\sigma}\ar[rr] && Y\ar[rr]^m &&E\ar[rr] && 0\\
 &&&&& E_\sigma\ar[ru]_{\iota} &&&&&
}$$
where $Y$ is the free $E\da E$-bimodule of rank $1$.  We check that the following diagram commutes:
$$\xymatrix@=34pt{
0\ar[r] &E\ar[r]^(0.4){\iota'}\ar[d]^{-1} &E\otimes_kE\ar[r]^{1\otimes x-x\otimes 1}\ar[d]^{\id\otimes\sigma} &E\otimes_k E\ar[r]^(0.6)m\ar@{=}[d] &E\ar[r]\ar@{=}[d] &0\\
0\ar[r] &E\ar[r]^(0.4){\iota_\sigma} &E\otimes_kE_\sigma\ar[r]^{\iota\circ m_\sigma} &E\otimes_k E\ar[r]^(0.6)m &E\ar[r] &0
}$$
where $\iota(e)=e\otimes x+ex\otimes 1$ and so our resolution is isomorphic to the resolution of $k[x]/\gen{x^{n+1}}$ constructed above for $n=1$.  Hence by Proposition \ref{comp} the periodic twist associated to the period $2$ resolution is naturally isomorphic to the square of the periodic twist associated to the period $1$ resolution.

When $n>1$, the kernel of the multiplication map $E\otimes_kE\to E$ is not even isomorphic to $E$ as a vector space, so we do not expect the analogue of the $\BP^{n}$-twist to decompose as the product of two autoequivalences.  However, by Theorem \ref{circles}, we do expect the functor to be given as the product of two combinatorial tilts.  We have a small circle of equivalences
$$\xymatrix{
\Db(A)\ar@/^3pc/[dd]^\sim \\
 \\
\Db(B)\ar@/^3pc/[uu]^\sim
}$$
where $B$ is in general not Morita equivalent to $A$.  When $n=1$ we have a similar circle of equivalences, but in this case $B\cong A$ and the derived equivalences are given by the analogues of spherical twists.

\subsection{Trivial extension algebras}
An infinite class of periodic symmetric algebras, indexed by simply laced Dynkin diagrams, was provided by Brenner, Butler, and King \cite{bbk}.  From this list, the only local algebra is $k[x]/\gen{x^2}$, which we have considered above, so this list of periodic algebras will only be applicable when we consider decomposable projective modules.

Let $A$ be the symmetric $k$-algebra generated by the paths of the quiver
$$\xymatrix@=10pt{
  Q_A&=&1\ar@/^/[rr]^{\alpha} & & 2\ar@/^/[ll]^\beta\ar@/^/[rr]^{\gamma} & & 3\ar@/^/[ll]^\delta
}$$
with relations $\alpha\beta\alpha=\beta\alpha\beta=\delta\gamma\delta=\gamma\delta\gamma=0$, $\alpha\gamma=\delta\beta=0$, and $\beta\alpha=\gamma\delta$.  $A$ is isomorphic to the trivial extension algebra of the path algebra of the Dynkin quiver of type $A_3$ with bipartitie orientation.  It is also isomorphic to the Brauer tree algebra of a line with 3 edges and no exceptional vertex, which is one of the algebras considered in \cite{rz}.

Each indecomposable projective $P_i$ associated to the vertex $i$ has endomorphism algebra $k[x]/\gen{x^2}$, and so we can define autoequivalences $\Psi_i=\Psi_{P_i}$ from the period $1$ twisted-periodic structure.   As described in both \cite{rz} and \cite{st}, these functors satisfy braid relations:
$$\Psi_1\Psi_2\Psi_1\cong \Psi_2\Psi_1\Psi_2\fs$$

Let $P=P_1\oplus P_2$ be the direct sum of the indecomposable projective modules associated to the first two vertices of the quiver of $A$.  Then $E=\End_A(P)^\op$ is isomorphic to the path algebra of the quiver
$$\xymatrix@=10pt{
  Q_E&=&1\ar@/^/[rr]^{\alpha} & & 2\ar@/^/[ll]^\beta
}$$
with relations $\alpha\beta\alpha=\beta\alpha\beta=0$.  $E$ is isomorphic to the trivial extension algebra of the path algebra of the Dynkin quiver of type $A_2$ and to the Brauer tree algebra with two edges and no exceptional vertex.

$E$ is twisted-periodic with periodicity $2$, with the automorphism $\sigma$ of $E$ induced by the graph automorphism of $Q_E$ which swaps the vertices $1$ and $2$.  This twisted-periodic algebra has truncated resolution
$$Q_{21}\oplus Q_{12}\arr{d} Q_{11}\oplus Q_{22}$$
where $Q_{ij}=Q_i\otimes_kQ_j^\vee$ and $Q_i$ denotes the projective $E$-module associated to the vertex $i$.  The map $d$ is given by 
$$(q\otimes r,s\otimes t)\mapsto (q\alpha\otimes e_1^\vee r + se_1\otimes \beta^\vee t, -qe_2\otimes\alpha^\vee r - s\beta\otimes e_2^\vee t)$$
where $q\otimes r\in Q_{21}$ and $s\otimes t\in Q_{12}$,
and so we obtain the two-sided tilting complex $X\in \Db(A\da A)$:
$$P_2\otimes_kP_1^\vee\oplus P_1\otimes_kP_2^\vee\arr{d} P_1\otimes_kP_1^\vee\oplus P_2\otimes_kP_2^\vee\arr{\ev}A$$
where the final map is the obvious evaluation map, and this induces the periodic twist $\Psi_P:\Db(A)\arr{\sim}\Db(A)$ which acts as $[2]$ on $P_1\oplus P_2$, sending $P_1$ to $P_2[2]$ and $P_2$ to $P_1[2]$, and acts as the identity on the simple module $P_3/\rad P_3$.

By direct calculation one can show that the functor $\Psi_P$ is naturally isomorphic to $\Psi_1\Psi_2\Psi_1$ (and hence also to $\Psi_2\Psi_1\Psi_2$) and so by using the projective module $P_1\oplus P_2$ we can realise the braid group action of the ``half-twist'' $s_1s_2s_1$ in a single step.

\subsection{Group algebras}
In order to construct new autoequivalences based on indecomposable projective modules $P$ we need to consider local twisted periodic algebras that are not isomorphic to $k[x]/\gen{x^n}$.

One way to obtain local algebras is to consider the group algebras $kG$ of finite $p$-groups $G$ over algebraically closed fields $k$ of characteristic $p$.  It is well-known that the unique simple module of such a group has a periodic projective resolution if and only if all abelian subgroups of $G$ are cyclic (see \cite[Section 4]{erdsko}), and by Theorem \ref{simpleper} this is equivalent to the group algebra $kG$ being twisted periodic.

Consider the group $G=Q_8$ of quaternions, a group of order $2^3$, over an algebraically closed field $k$ of characteristic $2$.  This group has a presentation
$$G=\gen{i,j|ij=j^{-1}i, ij=ji^{-1}}$$
and letting $E=kG$, we have an exact sequence of $A$-modules
$$0\to k\into E\delta\arr{d_3} E\gamma\oplus E\gamma' \arr{d_2} E\beta\oplus E\beta'\arr{d_1} E\alpha \to k \to0$$
which shows that the trivial $E$-module $k$ is periodic.  The first and last non-trivial maps are given by embedding the trivial module into the socle of $E\delta$ and by quotienting out by the radical of $E\alpha$.  Explicitly, the first maps $1\in k$ to $\sum_{g\in G}g$, and the last maps $g$ to $1\in k$ for all $g\in G$.  The differentials are given by
$$d_1(\beta)=(i+1)\alpha\textsf{; } d_1(\beta')=(j+1)\alpha\semic$$
$$d_2(\gamma)=(j+1)\beta + (ji+1)\beta'\textsf{; }d_2(\gamma')=(ij+1)\beta + (i+1)\beta'\semic$$
$$d_3(\delta)=(i+1)\gamma+(j+1)\gamma'\fs$$
See \cite{benper} for more details on this and similar examples, and \cite{erdsko} for more information on periodic groups.

Hence, given a symmetric $k$-algebra $A$ with a projective module $P$ such that $E=\End_A(P)^\op\cong kQ_8$, we can construct a derived autoequivalence $\Psi_P$ of $\Db(A)$ which acts as the shift $[4]$ on $P$ and which can be realised as the product of four combinatorial twists at $P$.

\end{document}